\input ./style/arxiv-general.cfg
\documentclass[aop,MSNbibl,dvips]{arximspdf}
\makeatletter
   \@ifpackageloaded{graphicx}{}{\usepackage{graphicx}}
\makeatother
\usepackage{mathrsfs}

\doi{10.1214/14-AOP942} 
\volume{43}
\issue{5}
\pubyear{2015}
\firstpage{2511}
\lastpage{2544}
\docsubty{FLA}

\makeatletter
\newcommand{\lleft}{\left}
\newcommand{\rright}{\right}
\newtheorem{theorem}{Theorem}[section]
\newtheorem{corollary}[theorem]{Corollary}
\newtheorem{thA}{Theorem}
\newtheorem{lemma}[theorem]{Lemma}
\newtheorem{proposition}[theorem]{Proposition}
\newproclaim{remark}[theorem]{Remark}

\newcommand{\binom}[2]{\pmatrix{#1 \cr #2}}
\newcommand{\underset}[2]{\mathop{#2}\limits_{#1}}
\newcommand{\pd}[2]{\frac{\partial#1}{\partial#2}}
\newcommand{\dd}{\partial}
\renewcommand{\dd}{\mathrm{d}}
\newcommand{\sign}{\operatorname{sign}}
\newcommand{\supp}{\operatorname{supp}}
\newcommand{\card}{\operatorname{card}}
\renewcommand{\P}{\mathsf{P}}
\newcommand{\E}{\mathsf{E}}
\newcommand{\Var}{\mathsf{Var}}
\renewcommand{\r}{\mathsf{r}}
\newcommand{\tA}{{\tilde{A}}}
\newcommand{\tB}{{\tilde{B}}}
\newcommand{\tc}{{\tilde{c}}}
\newcommand{\tPi}{{\tilde{\Pi}}}
\newcommand{\tH}{{\tilde{H}}}
\renewcommand{\Re}{\mathsf{Re}}
\renewcommand{\Re}{\operatorname{\mathfrak{Re}}}
\newcommand{\ZZ}{{\mathbf{Z}}}
\makeatother

\begin{document}
\begin{frontmatter}

\title{Exact Rosenthal-type bounds}
\runtitle{Exact Rosenthal-type bounds}

\begin{aug}
\author[A]{\fnms{Iosif}~\snm{Pinelis}\corref{}\ead[label=e1]{ipinelis@mtu.edu}}
\runauthor{I. Pinelis}
\affiliation{Michigan Technological University}
\address[A]{Department of Mathematical Sciences\\
Michigan Technological University\\
Houghton, Michigan 49931\\
USA\\
\printead{e1}}
\end{aug}
%

\received{\smonth{6} \syear{2013}}
\revised{\smonth{5} \syear{2014}}

%
\begin{abstract}
It is shown that, for any given $p\ge5$, $A>0$ and $B>0$, the exact
upper bound on
$\mathsf{E}|\sum X_i|^p$ over all independent zero-mean random variables
(r.v.'s) $X_1,\ldots,X_n$ such that $\sum\mathsf{E}  X_i^2=B$ and $\sum\mathsf{E}
|X_i|^p=A$ equals $c^p\mathsf{E} |\Pi_{\lambda}-\lambda|^p$, where $(\lambda
,c)\in(0,\infty)^2$ is the unique solution to the system of equations
$c^p\lambda=A$ and $c^2\lambda=B$, and $\Pi_\lambda$ is a Poisson
r.v. with mean $\lambda$.
In fact, a more general result is obtained, as well as other related ones.
As a tool used in the proof, a calculus of variations of 
moments of infinitely divisible distributions with respect to
variations of the L\'evy characteristics is developed.
\end{abstract}

%
\end{frontmatter}

\section{Introduction, summary, and discussion}\label{intro}
Let $\mathscr{X}$ denote the class of all finite sequences ${\mathbf{X}}
=(X_1,\ldots,X_n)$ of independent zero-mean random variables (r.v.'s).
For any ${\mathbf{X}}=(X_1,\ldots,X_n)\in\mathscr{X}$, let
%
\begin{equation}
\label{eqSX} S_{\mathbf{X}}:=X_1+\cdots+X_n.
\end{equation}
Take any real number
%
\begin{equation}
\label{eqp>2} p>2
\end{equation}
and any positive real numbers $A$ and $B$.
Consider
%
\begin{eqnarray}
\mathscr{X}_{p;A,B}&:=& \Biggl\{{\mathbf{X}}=(X_1,\ldots,X_n)\in\mathscr{X}\dvtx \sum_1^n
\E X_i^2=B, \sum_1^n
\E|X_i|^p=A \Biggr\}, \label{eqXX}
\\
\mathscr{X}_{p;\le A,\le B}&:=& \Biggl\{{\mathbf{X}}=(X_1,\ldots,X_n)\in\mathscr{X}\dvtx \sum_1^n
\E X_i^2\le B, \sum_1^n
\E|X_i|^p\le A \Biggr\}, \label{eqXX<}
\\
\mathscr{E}_{p;A,B}&:=&\sup\bigl\{\E|S_{\mathbf{X}}|^p
\dvtx {\mathbf{X}}\in\mathscr{X}_{p;A,B}\bigr\}, \label{eqEEA,B}
\\
\qquad \mathscr{E}_{p;\le A,\le B}&:=&\sup\bigl\{\E|S_{\mathbf{X}}|^p\dvtx
{\mathbf{X}}\in\mathscr{X}_{p;\le A,\le B}\bigr\}. \label{eqEE<A,<B}
\end{eqnarray}

Rosenthal's upper bound (Theorem~3 in \cite{rosenthal}) can be presented
by the inequality
%
\begin{equation}
\label{eqros} \mathscr{E}_{p;A,B}\le C_p \max
\bigl(A,B^{p/2}\bigr),
\end{equation}
with $C_p:=(p/2)^{p/2} 2^{p+p^2/4}$. In particular, this implies that
$\mathscr{E}_{p;A,B}<\infty$. For
some of the subsequent developments, see, for example, Sections~4 and 5 in \cite{pin94}, \cite{sibam,ibr-sankhya} and
references therein.

%
\begin{proposition}\label{propEE}
One has $\varnothing\ne\mathscr{X}_{p;A,B}\subseteq\mathscr{X}_{p;\le A,\le B}$.
Moreover, one has the homogeneity property
$\mathscr{E}_{p;\kappa^p A,\kappa^2 B}
=\kappa^p\mathscr{E}_{p;A,B}$ for all real $\kappa>0$.
Furthermore,
$\mathscr{E}_{p;A,B}$ is nondecreasing in $A$ and in $B$ and hence
%
\begin{equation}
\label{eqEEle=} \mathscr{E}_{p;\le A,\le B}=\mathscr{E}_{p;A,B}.
\end{equation}
\end{proposition}

All the necessary proofs are deferred to Sections~\ref{proofs} and~\ref{proofs,props}.
In particular, Proposition~\ref{propEE} will be proved in Section~\ref{proofs,props}.

Using
Proposition~\ref{propEE}, one can easily see
(cf.~\cite{ibr-sankhya})
that the problem of finding a good expression of $\mathscr{E}_{p;A,B}$
is equivalent to that of finding, for an arbitrary balancing parameter
$\gamma\in(0,\infty)$, a good expression of the best constant
$C_{p;\gamma}$ in the Rosenthal-type inequality
%
\begin{equation}
\label{eqga} \mathscr{E}_{p;A,B}\le C_{p;\gamma}\max\bigl(\gamma
A,B^{p/2}\bigr);
\end{equation}
cf. (\ref{eqros}).
Indeed, one has:

%
\begin{proposition}\label{propEE,C}
$
C_{p;\gamma}=\mathscr{E}_{p;1/\gamma,1} \mbox{ and } \mathscr{E}_{p;A,B}=B^{p/2}C_{p;B^{p/2}/A}$.
\end{proposition}

The idea of balancing the contributions of the terms $A$ and $B^{p/2}$
in the Rosen\-thal-type bounds, depending on the relative sizes of these
terms, goes back at least to the Corollary in \cite{pin80}; see also
Sections~4 and 5, Remark~6.8, and Theorem~8.3 in \cite{pin94},
Proposition~9.2 in \cite{pin94-arxiv},
Corollaries~3.1, 3.2 in  \cite{re-centerpubl}, and
Corollaries~2,~3,~4 in \cite{pin12-2smooth}.

For any real $\lambda>0$, let $\Pi_\lambda$ denote a r.v. with the
Poisson distribution with mean~$\lambda$, and then introduce the
corresponding centered r.v.
\[
\tPi_\lambda:=\Pi_\lambda-\lambda.
\]

Using Theorem~4 by Utev \cite{utev-extr}, Bestsennaya and Utev \cite
{bestsen} showed that
%
\begin{equation}
\label{eqrosenthal-conjec} \mathscr{E}_{p;A,B}= c^p\E|
\tPi_\lambda|^p\qquad\mbox{if }p=4,6,\ldots,
\end{equation}
where
%
\begin{equation}
\label{eqla,c} \lambda:=\lambda_p(A,B):= \biggl(\frac{B^{p/2}}A
\biggr)^{2/(p-2)}\quad\mbox{and}\quad c:=c_p(A,B):= \biggl(
\frac{A}{B} \biggr)^{1/(p-2)},\hspace*{-30pt}
\end{equation}
so that the pair $(\lambda,c)\in(0,\infty)^2$ is the unique solution
to the system of equations
\[
c^2\lambda=B\quad\mbox{and}\quad c^p\lambda=A.
\]
Obviously, if $p$ is an even natural number, then the absolute $p$th
moment $\E|X|^p$ of a r.v. $X$ is the same as its $p$th moment $\E
X^p$. This fact allows the proof in \cite{bestsen} to be based on the
well-known representation of moments in terms of cumulants and the
log-convexity of $\int_\mathbb{R}|x|^r G(\dd x)$ in $r>0$, for any
nonnegative measure $G$.

Under the additional restriction that the $X_i$'s be symmetric(ally
distributed), exact Rosenthal-type
bounds were obtained in \cite{utev-extr,zinnetal,ibr-shar97,ibr-sankhya}.
In particular, it was shown by Utev \cite{utev-extr} that
%
\begin{equation}
\label{equtev,p} \mathop{\sup_{{\mathbf{X}}\in\mathscr{X}_{p;\le A,\le
B},}}_{{\mathbf{X}}\ \mathrm
{is\ symmetric}}
\E|S_{\mathbf{X}}|^p =\mathop{\sup_{{\mathbf{X}}\in\mathscr{X}_{p;A,B},}}_{{\mathbf
{X}}\ \mathrm{is\
symmetric}}
\E|S_{\mathbf{X}}|^p =c^p\E\bigl|\Pi_{\lambda/2}-
\Pi_{\lambda/2}^\diamond\bigr|^p
\end{equation}
if $p>4$,
where $\lambda$ and $c$ are as in (\ref{eqla,c}), and $\Pi_{\lambda
/2}^\diamond$ is an independent copy of $\Pi_{\lambda/2}$.\vspace*{6pt}


Take
any
%
\begin{equation}
\label{eqq} q\in(2,p] 
\end{equation}
and then take any r.v. $X$ such that
%
\begin{equation}
\label{eqX} \E|X|^q<\infty.
\end{equation}
Consider
%
\begin{eqnarray}
\label{eqXXX} \mathscr{X}_{p;X;A,B}&:=& \{{\mathbf{X}}\in\mathscr{X}_{p;A,B}\dvtx {\mathbf{X}} \mbox{ is independent of }X \},
\nonumber\\[-8pt]\\[-8pt]
\mathscr{X}_{p;X;\le A,\le B}&:=& \{{\mathbf{X}}\in\mathscr{X}_{p;\le
A,\le B}\dvtx {\mathbf{X}}\mbox{ is independent of }X \}.\nonumber
\end{eqnarray}

The main result of the present paper is:

%
\begin{theorem}\label{th}
Suppose that $p\ge q\ge5$ and $\E X=0$.
Then
%
\begin{eqnarray}\label{eq}
\sup_{{\mathbf{X}}\in\mathscr{X}_{p;X;\le A,\le B}}\E
|X+S_{\mathbf{X}}|^q
&=& \sup_{{\mathbf{X}}\in\mathscr{X}_{p;X;A,B}}\E|X+S_{\mathbf{X}}|^q
\nonumber\\[-8pt]\\[-8pt]
&=&\max\bigl(\E|X+c\tPi_\lambda|^q, \E|X-c
\tPi_\lambda|^q \bigr),\nonumber
\end{eqnarray}
where $\lambda$ and $c$ are as in (\ref{eqla,c}), and the r.v. $\tPi
_\lambda$ is independent of $X$.
\end{theorem}


In the special case when $X=0$ and $q=p$, Theorem~\ref{th} yields
\[
\mathscr{E}_{p;A,B}=\mathscr{E}_{p;\le A,\le B}= c^p\E|
\tPi_\lambda|^p\qquad\mbox{if }p\ge5;
\]
cf. (\ref{eqrosenthal-conjec}).

Allowing $q$ in Theorem~\ref{th} to differ from $p$ not only provides
a more general result, but also helps with the proof. Indeed,
Theorem~\ref{th} will be first proved in the case when $p>q>5$ [see
(\ref{eq5<q<p})], and then the proof will be completed by limit
transitions in $q$ and in $p$.

%
\begin{remark}\label{remgeneralf}
It is of substantial interest to obtain exact Rosenthal-type
inequalities for moment functions more general than the function
$|\cdot|^p$ used in Theorem~\ref{th}; cf., for example, \cite{burk,zinnetal}.
In fact, one can indeed easily extend the result of Theorem~\ref{th}
to the class of all moment functions of the form
%
\begin{eqnarray}
\label{eqmf} x&\longmapsto&\int_{[5,p]\times[0,\infty)}(a+x)_+^r
\nu_1(\dd r\times\dd a)
\nonumber\\[-8pt]\\[-8pt]
&&{} +\int_{[5,p]\times[0,\infty)}(a-x)_+^r
\nu_2(\dd r\times\dd a),\nonumber
\end{eqnarray}
where $\nu_1$ and $\nu_2$ are any nonnegative Borel measures on the
set $[5,p]\times[0,\infty)$ such that the resulting moment function
is real-valued; of course, the moment function $x\mapsto|x|^p
(=x_+^p+(-x)_+^p )$
is just one member of this class;
as usual, we let $x_+:=0\vee x$ and $x_+^r:=(x_+)^r$ for all real $x$
and all real $r>0$.
To see why this extension of Theorem~\ref{th} is valid, one needs to
look at the place in the proof of the theorem that imposes the
narrowest restriction on the moment function---which is the condition
that the difference $h''(u\alpha s)-u^{p-4}h''(\alpha s)$, considered
in (\ref{eq>0}), be strictly positive for all $u$, $\alpha$, and $s$
in $(0,1)$.
The class of functions given by (\ref{eqmf}) may be compared with
classes of moment functions considered, for example, in \cite{T2,pin98,asymm}.
\end{remark}

In what follows, to avoid repetitiveness, it is assumed that the
different instances of
\[
\framebox{\mbox{all r.v.'s entering the same expression are
independent.}}
\]
Thus, conditions such as that of the independence of the r.v.'s $\tPi
_\lambda$ and $X$ in Theorem~\ref{th} may not be explicitly stated in
the sequel.

Theorem~\ref{th} is complemented by:

%
\begin{theorem}\label{th2<p<3}
Suppose that $p\in(2,3]$ and $\E|X|^p<\infty$ (the condition \mbox{$\E
X=0$} is not needed here).
Then
%
\begin{eqnarray}\label{eq2<p<3}
\sup_{{\mathbf{X}}\in\mathscr{X}_{p;X;\le A,\le B}}\E
|X+S_{\mathbf{X}}|^p
&=& \sup_{{\mathbf{X}}\in\mathscr{X}_{p;X;A,B}}\E|X+S_{\mathbf{X}}|^p
\nonumber\\[-8pt]\\[-8pt]
&=&A+\E\bigl|X+B^{1/2}Z\bigr|^p.\nonumber
\end{eqnarray}
Here and in what follows, $Z\sim N(0,1)$, unless specified otherwise.
\end{theorem}

Theorem~\ref{th2<p<3} is based on a result by Tyurin \cite{tyurinSPL}.
In the case $p=3$, important for applications to Berry--Esseen bounds,
a certain refinement of (\ref{eq2<p<3}) was obtained in
Corollary~2 from \cite{p=3publ}, based on the main result in the paper \cite
{pin-hoeff-arxiv}, a shorter version of which appeared in \cite
{pin-hoeff-published}.


One has the following interpretation of the last expression in (\ref
{eq2<p<3}), in terms of centered Poisson r.v.'s $\tPi_{\lambda_1}$
and $\tPi_{\lambda_2}$ (such that the r.v.'s $X$, $\tPi_{\lambda
_1}$, and $\tPi_{\lambda_2}$ are independent).

%
\begin{proposition}\label{propA+E}
Suppose that $p\in(2,3]$ and $\E|X|^p<\infty$. Then
%
\begin{eqnarray}\label{eq=lim}
&& A+\E\bigl|X+B^{1/2}Z\bigr|^p\nonumber
\\
&&\qquad =\lim\bigl(\E|X+c_1\tPi_{\lambda_1}+c_2
\tPi_{\lambda_2}|^p\dvtx
\\
&&\hspace*{52pt}  (c_1,c_2,\lambda_1,\lambda_2)\in Q_{p;A,B}, c_1\to0, |c_2|\to\infty\bigr),\nonumber
\end{eqnarray}
where
%
\begin{eqnarray}\label{eqQ=}
Q_{p;A,B} &:=& \bigl\{(c_1,c_2,
\lambda_1,\lambda_2)\in\mathbb{R}^2\times(0,
\infty)^2\dvtx
\nonumber\\[-8pt]\\[-8pt]
&&\hspace*{5pt}  c_1^2\lambda_1+c_2^2 \lambda_2=B, |c_1|^p\lambda_1+|c_2|^p\lambda_2=A \bigr\}.\nonumber
\end{eqnarray}
\end{proposition}

Proposition~\ref{propA+E} will be useful in the proof of Theorem~\ref{th2<p<3}.

Now one can present a unified form
of
the exact upper bounds in (\ref{eq}) and (\ref{eq2<p<3}):

%
\begin{corollary}\label{corA+E}
Suppose that $p\in(2,3]\cup[5,\infty)$ and $\E|X|^p<\infty$. For
$p\in[5,\infty)$, also suppose that $\E X=0$. Then
%
\begin{eqnarray}
\label{equnified}
&& \sup_{{\mathbf{X}}\in\mathscr{X}_{p;X;\le A,\le B}}\E
|X+S_{\mathbf{X}}|^p\nonumber
\\
&&\qquad = \sup_{{\mathbf{X}}\in\mathscr{X}_{p;X;A,B}}\E|X+S_{\mathbf{X}}|^p
\\
&&\qquad =\sup\bigl\{\E|X+c_1\tPi_{\lambda_1}+c_2
\tPi_{\lambda
_2}|^p\dvtx (c_1,c_2,
\lambda_1,\lambda_2)\in Q_{p;A,B} \bigr\}.\nonumber
\end{eqnarray}
%
\end{corollary}

By Theorem~\ref{th2<p<3} and Proposition~\ref{propA+E}, for $p\in
(2,3]$ the last supremum in (\ref{equnified}) is ``attained in the
limit'' as $c_1\to0$ and $ |c_2|\to\infty$, whereas, by Theorem~\ref{th}, for $p\ge5$ the same supremum is (actually) attained at
$(c_1,c_2,\lambda_1,\lambda_2)=(c,0,\lambda,0)$ or at
$(c_1,c_2,\lambda_1,\lambda_2)=(-c,0,\lambda,0)$, where $\lambda$
and $c$ are as in (\ref{eqla,c}).

The cases $p\in(3,4)$ and $p\in(4,5)$ remain open. Certain
considerations suggest that Theorem~\ref{th2<p<3} should hold for
$p\in(3,4)$ as well, whereas Theorem~\ref{th} should hold for $p\in
(4,5)$---at least when $q=p$.
For $q=p=4$, it is easy to see that the ``answers'' in (\ref{eq}) and
(\ref{eq2<p<3}) coincide with each other:
\[
\E|X+c\tPi_\lambda|^4=\E|X-c\tPi_\lambda|^4=A+
\E\bigl|X+B^{1/2}Z\bigr|^4.
\]

This situation may be compared with the one concerning the exact\break
Khinchin-type upper bound. There the summands are weighted independent
Rademacher r.v.'s $X_1=a_1\varepsilon_1,\ldots,X_n=a_n\varepsilon_n$,
where $\P(\varepsilon_i=\pm1)=1/2$, and the weights $a_1,\ldots,a_n$
are real numbers subject to the restriction $\sum_1^n a_i^2=1$.
Since these summands have each a simplest symmetric distribution, and
there is only one restriction here on the sum of the moments, $\sum
_1^n\E X_i^2=\sum_1^n a_i^2$, it appears that the problem of the exact
Khinchin-type upper bound is significantly simpler than its
Rosenthal-type counterpart. Indeed, in 1960 Whittle \cite{whittle}
gave a very simple proof of the exact Khinchin-type upper bound, $\E
|Z|^p$, for the case $p\ge3$. The proof in \cite{whittle} was based
on the fact that, again for $p\ge3$, the second derivative of $|x|^p$
in $x$ is convex in $x\in\mathbb{R}$. It was claimed in \cite
{whittle} that the result holds for all real $p\ge2$, but that was not
supported by the proof.
Actually, the problem of the exact Khinchin-type upper bound in the
case $p\in(2,3)$ turned to be very difficult and was solved only in
1981 by Haagerup \cite{haag}. Haagerup's proof was somewhat simplified
in \cite{naz-podk}; see also \cite{dichot}.
One may speculate that the case $p\ge5$ in the Rosenthal-type context
is parallel to the case $p\ge3$ in the Khinchin-type one, whereas the
Rosenthal-type case of a small noninteger $p\in(3,4)\cup(4,5)$ is
parallel to the Khinchin-type case of $p\in(2,3)$.
If so, the remaining Rosenthal-type case of $p\in(3,4)\cup(4,5)$ may
be exceedingly difficult, on comparing the treatment of the
Rosenthal-type case of $p\ge5$ in the present paper with that of the
Khinchin-type case of $p\ge3$ in \cite{whittle}.
One may also note here that the condition $p\ge5$ will be used twice,
and in rather different ways, in the proof of Theorem~\ref{th}, namely
in the proofs of Propositions~\ref{prop<2} and~\ref{propcard=1}.

For the symmetric case, one has:

%
\begin{theorem}\label{thsymm}
Suppose that $p\ge q\ge5$ and $\E X=0$.
Then
%
\begin{eqnarray}\label{eqsymm}
&& \mathop{\sup_{{\mathbf{X}}\in\mathscr{X}_{p;X;\le A,\le B},}}_{{\mathbf{X}}\ \mathrm
{is\ symmetric}}
\E|X+S_{\mathbf{X}}|^q
\nonumber\\[-8pt]\\[-8pt]
&&\qquad =\mathop{\sup_{{\mathbf{X}}\in\mathscr{X}_{p;X;A,B},}}_{{\mathbf
{X}}\ \mathrm{is\ symmetric}}
\E|X+S_{\mathbf{X}}|^q =\E\bigl|X+c\Pi_{\lambda/2}-c
\Pi_{\lambda/2}^\diamond\bigr|^q,\nonumber
\end{eqnarray}
where $\lambda$ and $c$ are as in (\ref{eqla,c}) and, as in (\ref
{equtev,p}), $\Pi_{\lambda/2}^\diamond$ is an independent copy of
$\Pi_{\lambda/2}$.
\end{theorem}

Theorem~\ref{thsymm} generalizes (\ref{equtev,p}), but only for $p\ge5$.
The generalization has two aspects: (i) letting $q$ differ from $p$ and
(ii) introducing the extra summand $X$. Note that $X$ is not required
to be symmetric in Theorem~\ref{thsymm}.

An advantage of having the extra summand $X$ is illustrated by the
following straightforward combination of Theorems~\ref{th} and \ref{thsymm}.

%
\begin{corollary}\label{cor}
Suppose that $p\ge q\ge5$ and $\E X=0$.
Take any positive real numbers $A_0,B_0,A_1,B_1$.
For each $j\in\{0,1\}$, let $\lambda_j:=\lambda_p(A_j,B_j)$ and
$c_j:=c_p(A_j,B_j)$, in accordance with (\ref{eqla,c}).
Then
%
\begin{eqnarray}\label{eqcomb}
\mathop{\mathop{\mathop{\sup_{{\mathbf{X}}\in
\mathscr{X}_{p;\le
A_0,\le B_0},}}_{{\mathbf{Y}}\in\mathscr{X}_{p;\le A_1,\le B_1},}}_
{{\mathbf{X}}\ \mathrm{is\ symmetric},}}_{X,{\mathbf{X}},{\mathbf
{Y}}\ \mathrm{are\ independent}}
\E|X+S_{\mathbf{X}}+S_{\mathbf{Y}}|^q
&=& \mathop{\mathop{\mathop{
\sup_{{\mathbf{X}}\in\mathscr{X}_{p;A_0,B_0},}}_{{\mathbf{Y}}\in\mathscr{X}_{p;A_1,B_1},}}_{{\mathbf{X}}\ \mathrm
{is\ symmetric},}}_{X,{\mathbf{X}},{\mathbf{Y}}\ \mathrm{are\
independent}}
\E|X+S_{\mathbf{X}} +S_{\mathbf{Y}}|^q\hspace*{-8pt}\nonumber
\\
&=& \E\bigl|X+c_0\Pi_{\lambda_0/2}-c_0
\Pi_{\lambda_0/2}^\diamond+c_1\tPi_{\lambda_1}
\bigr|^q \hspace*{-8pt}
\\
&&{} \vee\E\bigl|X+c_0\Pi_{\lambda_0/2}-c_0
\Pi_{\lambda_0/2}^\diamond-c_1\tPi_{\lambda_1}
\bigr|^q.\hspace*{-8pt}\nonumber
\end{eqnarray}
\end{corollary}

This follows immediately from Theorems~\ref{thsymm} and \ref{th}, by
taking first the supremum in ${\mathbf{X}}$ (say) and then in
${\mathbf{Y}}$.

Note that, in the case when $X$ is symmetric (or, in particular, zero),
the maximum in (\ref{eqcomb}) simplifies to $\E|X+c_0\Pi
_{\lambda_0/2}-c_0\Pi_{\lambda_0/2}^\diamond+c_1\tPi_{\lambda
_1} |^q$.

Corollary~\ref{cor} may be useful when some, but not all, of the
independent summands are known to be symmetric.

For the calculation of absolute moments, especially such more
complicated ones as in the maximum expression in (\ref{eqcomb}),
Fourier- or Fourier--Laplace-type identities such as those given in
\cite{positive} can be effective; one of such identities will be
reproduced in the present paper as (\ref{eqstart}).

\section{Proof of Theorem~\texorpdfstring{\protect\ref{th}}{1.3}}\label{proofs}

\subsection{Domination by the accompanying compound Poisson distribution}\label{accomp}

\renewcommand{\thethA}{A}
\begin{thA}\label{teoA}
Let $f\dvtx \mathbb{R}\to\mathbb{R}$ be any
twice continuously differentiable function such that $f$ and $f''$ are convex.
Let $G$ be any finite nonnegative Borel measure on $\mathbb{R}$ such
that $G(\{0\})=0$ and $\int_\mathbb{R}xG(\dd x)=0$, and then let
$X_G$ be any r.v. with the characteristic function $t\mapsto\exp\int
_\mathbb{R}(e^{itx}-1)G(\dd x)$.
Then
\[
\label{equtev} \sup\bigl\{\E f(S_{\mathbf{X}})\dvtx {\mathbf{X}}\in
\mathscr{X}, G_{\mathbf{X}}=G\bigr\}=\E f(X_G),
\]
where $S_{\mathbf{X}}$ is as in (\ref{eqSX}) and $G_{\mathbf{X}}$ is
the ``sum of the
tails'' measure defined by
\[
G_{\mathbf{X}}(E):=\sum\P\bigl(X_i\in E\setminus\{0\}
\bigr)
\]
for all Borel subsets $E$ of $\mathbb{R}$.
In particular, for all $x\in\mathbb{R}$ and all real $p\ge3$,
\begin{eqnarray*}
\sup\bigl\{\E|S_{\mathbf{X}}-x|^p\dvtx {\mathbf{X}}\in\mathscr{X},
G_{\mathbf{X}}=G\bigr\}&=&\E|X_G-x|^p,
\\
\sup\bigl\{\E(S_{\mathbf{X}}-x)_+^p\dvtx {\mathbf{X}}\in
\mathscr{X}, G_{\mathbf{X}}=G\bigr\}&=&\E(X_G-x)_+^p.
\end{eqnarray*}
\end{thA}

Theorem~\ref{teoA} is essentially the same as the mentioned Theorem~4 by Utev
\cite{utev-extr}; cf. \cite{prokhorov,pin-utev84,spher}. The
assumptions on $f$ in Theorem~4 from \cite{utev-extr}, were slightly
different; namely, it was assumed there that $f''$ is convex whereas
$f$ is nonnegative and satisfies a certain limited growth condition,
which latter may be dropped, by Proposition~1 and
Lemma~4 in \cite{asymm}, provided that $f$ and $f''$ are convex, as in Theorem~\ref{teoA}.

\begin{remark*}
If a r.v. $X$ has a finite expectation and a function $f\dvtx \mathbb
{R}\to\mathbb{R}$ is convex, then, by Jensen's inequality, $\E f(X)$
always exists in $(-\infty,\infty]$.
\end{remark*}

Let us complement Theorem~\ref{teoA} by the following standard lemma;
cf., for example, \cite{pin-utev84,pin-utev-exp} or the paragraphs
containing formulas (6.1) and (6.2) in \cite{pin94}.

%
\begin{lemma}\label{lemlim}
Let $G$ be any finite nonnegative Borel measure on $\mathbb{R}$ such that
$\int_\mathbb{R}|x|^p G(\dd x)=A$ and
$\int_\mathbb{R}x^2 G(\dd x)=B$; such a measure $G$ exists.
Let then $X_G$ be any r.v. with the characteristic function $t\mapsto
\exp\int_\mathbb{R}(e^{itx}-1-itx)G(\dd x)$.
Then there exists a sequence $(\ZZ_n)$ in $\mathscr{X}_{p;A,B}$ such
that $S_{\ZZ_n}\stackrel{\mathrm{D}}\longrightarrow X_G$, where
$\stackrel{\mathrm{D}}\longrightarrow$ denotes the convergence in
distribution.
In particular, it follows that $\mathscr{X}_{p;A,B}\ne\varnothing$.
\end{lemma}

The conditions on $G$ in Lemma~\ref{lemlim} are different from those
in Theorem~\ref{teoA}. In particular, the conditions
$G(\{0\})=0$ and $\int_\mathbb{R}xG(\dd x)=0$ are not required in
Lemma~\ref{lemlim}.
However, when the condition $\int_\mathbb{R}xG(\dd x)=0$ does hold,
the definition of the r.v. $X_G$ in Lemma~\ref{lemlim} is consistent
with that in Theorem~\ref{teoA}.
Also, the condition $p\ge5$ imposed in Theorem~\ref{th} is not needed
in Lemma~\ref{lemlim}; rather, it is enough to assume there that the
general condition (\ref{eqp>2}) holds.

\begin{pf*}{Proof of Lemma~\ref{lemlim}}
First, concerning the existence of $G$, note that all the conditions on
$G$ imposed in Lemma~\ref{lemlim} are satisfied by the measure
$\lambda\delta_c$, where $\lambda$ and $c$ are as in (\ref{eqla,c})
and $\delta_u$ denotes the Dirac probability measure at $u$.

Next, for each natural $n$ and all $j\in\{1,\ldots,n\}$, let
%
\begin{equation}
\label{eqZ} Z_{j,n}:=W_{j,n}-\E W_{j,n},
\end{equation}
where the $W_{j,n}$'s are independent identically distributed r.v.'s
with the distribution determined by the condition that
%
\begin{equation}
\label{eqWdistr} \E f(W_{j,n})=f(0)+\frac{\kappa_n}n \int
_\mathbb{R}\bigl[f(\gamma_n x)-f(0)\bigr]G(\dd x)
\end{equation}
for all (say) bounded or nonnegative Borel functions $f\dvtx \mathbb
{R}\to\mathbb{R}$, and where in turn $\kappa_n$ and $\gamma_n$ are
positive real numbers such that
$\frac{\kappa_n}n \int_\mathbb{R}G(\dd x)\le1$;
the latter condition
is precisely what is needed for formula (\ref{eqWdistr}) to define a
probability distribution.
It follows that for $r\in\{2,p\}$,
\begin{eqnarray*}
&& \sum_1^n\E|Z_{j,n}|^r=n
\E|Z_{1,n}|^r=F_r \biggl(\frac1n,\kappa
_n,\gamma_n \biggr),
\end{eqnarray*}
where
\[
F_r(\alpha,\kappa,\gamma):= |\kappa\gamma
m_G|^r|\alpha|^{r-1}\sign\alpha+\kappa
\gamma^r\int_\mathbb{R} \bigl(|x-\kappa\alpha
m_G|^r-|\kappa\alpha m_G|^r
\bigr) G(\dd x) %
\]
%
and $m_G:=\int_\mathbb{R}x G(\dd x)$.
Introducing now the vector function ${\mathbf{F}}:=(F_2,F_p)$, we see
that it is\vspace*{4pt} continuously differentiable on $\mathbb{R}\times(0,\infty
)^2$, and the Jacobian matrix
$
{\fontsize{7pt}{10pt}\selectfont{\pmatrix{
\hspace*{-2pt}\!\pd{F_2}{\kappa}\! & \pd{F_2}{\gamma}\hspace*{-2pt} \cr
\hspace*{-2pt}\pd{F_p}\kappa\! & \pd{F_p}\gamma\hspace*{-2pt}}}}$ 
at the point $(\alpha,\kappa,\gamma)=(0,1,1)$ is
$
{\fontsize{8.2pt}{10pt}\selectfont{\pmatrix{
\hspace*{-2pt}B\! &\! 2B\hspace*{-2pt} \cr
\hspace*{-2pt}A\! &\! pA\hspace*{-2pt}}}}$,
which is\vspace*{4pt} nonsingular.
Moreover,
${\mathbf{F}}(0,1,1)=(B,A)$.
So, by the implicit function theorem, there exist a positive real
number $\alpha_0$ and continuously differentiable functions $\tilde
\kappa\dvtx (-\alpha_0,\alpha_0)\to\mathbb{R}$ and $\tilde\gamma
\dvtx (-\alpha_0,\alpha_0)\to\mathbb{R}$ such that $\tilde\kappa
(0)=\tilde\gamma(0)=1$ and ${\mathbf{F}} (\alpha,\tilde\kappa
(\alpha),\tilde\gamma(\alpha) )=(B,A)$ for all $\alpha\in
(-\alpha_0,\alpha_0)$.
For all
natural $n>1/\alpha_0$, letting now
$
\kappa_n:=\tilde\kappa(\frac1n)$ and $\gamma_n:=\tilde\gamma
(\frac
1n)$, 
one sees that
$\sum_1^n\E|Z_{j,n}|^2=B$ and $\sum_1^n\E|Z_{j,n}|^p=A$, so that
$
\ZZ_n:=(Z_{1,n},\ldots,Z_{n,n})\in\mathscr{X}_{p;A,B}$. %
Thus, indeed $\mathscr{X}_{p;A,B}\ne\varnothing$.

Moreover, $\kappa_n\to\tilde\kappa(0)=1$ and $\gamma_n\to\tilde
\gamma
(0)=1$ (the convergence in this context is of course as $n\to\infty$).
So, by (\ref{eqZ}) and (\ref{eqWdistr}),
\begin{eqnarray*}
&& \E\exp(itS_{\ZZ_n} ) = \biggl[1+\frac{\kappa_n}n \int
_\mathbb{R} \bigl(e^{it\gamma_n
x}-1 \bigr) G(\dd x)
\biggr]^n e^{-it\kappa_n\gamma_n m_G}
\\
&&\quad \longrightarrow\quad \exp\int_\mathbb{R}\bigl(e^{itx}-1-itx
\bigr)G(\dd x) =\E\exp(itX_G )
\end{eqnarray*}
for all real $t$, so that indeed $S_{\ZZ_n}\stackrel{\mathrm
{D}}\longrightarrow X_G$.
\end{pf*}

\subsection{Zero-mean truncation of zero-mean r.v.'s}\label{trunc}

%
\begin{proposition}\label{proptrunc}
Let $Y$ be any zero-mean r.v. Then for any real $M>0$ there is an r.v.
$Y_M$ with the following properties:
\begin{longlist}[(iii)]
\item[(i)] $\E Y_M=0$;
\item[(ii)] $|Y_M|\le M\wedge|Y|$;
\item[(iii)] $\E f(Y_M)\le\E f(Y)$ for all convex functions $f\dvtx \mathbb
{R}\to\mathbb{R}$;
\item[(iv)] $Y_M\to Y$ almost surely (a.s.) as $M\to\infty$.
\end{longlist}
\end{proposition}

This follows immediately from Proposition~3.15 in \cite{disintegr}, and
Jensen's inequality
on letting
\[
Y_M:=Y \operatorname{\mathrm{I}}(E_M)=\E(Y|
\mathscr{F}_M),
\]
where $E_M:=\{|Y|\le M,|\r(Y,U)|\le M\}$,
$U$ is any r.v. which is independent of $Y$ and uniformly distributed
on the unit interval
$[0,1]$,
$\r$ stands for
the reciprocating function of (the distribution of) the r.v. $Y$ in
accordance with the definition (formula (2.6) in \cite{disintegr}), and $\mathscr{F}_M$ is the $\sigma$-algebra generated by all events of the form
$E_M\cap\{Y\le y,U\le u\}$ with any real $y$ and $u$.
Note that, by Proposition~3.6 in \cite{disintegr}, $|\r(Y,U)|<\infty$ a.s.
The r.v. $U$, which may be referred to as a randomizing r.v., is used
to split atoms of the distribution of $Y$, as such splitting may be
needed to satisfy the condition $\E Y_M=0$.

\subsection{Differentiation under the integral sign}\label{diff}

Take any measurable space $(\Omega,\mathscr{F})$ with a measure $\mu
\dvtx \mathscr{F}\to\mathbb{C}$.
Take also any $t_*\in(0,\infty)$.
Let $f\dvtx \Omega\times[0,t_*)\to\mathbb{R}$.
Suppose that for each $t\in[0,t_*)$ the function $\Omega\ni\omega
\mapsto f(\omega,t)$ is $\mu$-integrable, and let
\[
F(t):=\int_\Omega\mu(\dd\omega) f(\omega,t).
\]
Suppose also that, for each $\omega\in\Omega$, the function
$[0,t_*)\ni t\mapsto f(\omega,t)$ is continuous and has a
right-continuous right-hand side derivative $[0,t_*)\ni t\mapsto
(\partial_2f) (\omega,t)\in[-\infty,\infty]$
such that the function $\Omega\ni\omega\mapsto(\partial_2f)
(\omega,t)$ is $\mathscr{F}$-measurable, for each $t\in[0,t_*)$.

%
\begin{lemma}\label{lemunderint}
Suppose that for each pair $(t,\varepsilon)\in[0,t_*)\times(0,\infty
)$ there exist a set $\Omega_{t,\varepsilon}\in\mathscr{F}$, a
measurable
function $g_{t,\varepsilon}\dvtx \Omega_{t,\varepsilon}\to
[0,\infty]$ and a real number $h_{t,\varepsilon}\in(0,t_*-t)$ such
that $\int_{\Omega_{t,\varepsilon}}|\dd\mu| g_{t,\varepsilon
}<\infty$,
%
\begin{equation}
\label{eq<g} \bigl|(\partial_2f) (\omega,v) \bigr|\le g_{t,\varepsilon}(
\omega)\qquad\mbox{for all } (\omega,v)\in\Omega_{t,\varepsilon}
\times[t,t+h_{t,\varepsilon})
\end{equation}
and
%
\begin{equation}
\label{eqsupto0} \sup_{v\in[t,t+h_{t,\varepsilon})}\int_{\Omega\setminus
\Omega
_{t,\varepsilon}} \bigl|\mu(
\dd\omega) (\partial_2f) (\omega,v) \bigr| \underset{\varepsilon
\downarrow0} {\longrightarrow}0.
\end{equation}
Then
%
\begin{equation}
\label{eqF=30} \qquad F'(t+):=\lim_{h\downarrow0}
\frac{F(t+h)-F(t)}h =I(t):=\int_\Omega\mu(\dd\omega) (
\partial_2f) (\omega,t)\in\mathbb{C}.
\end{equation}
\end{lemma}

The following lemma is
a special case of Lemma~\ref{lemunderint}.

%
\begin{lemma}\label{lemunderint,simpl}
Suppose that there exists a $\mu$-integrable function $g\dvtx \Omega
\to[0,\infty]$ such that
%
\begin{equation}
\label{eq<g,simpl} \bigl|(\partial_2f) (\omega,t) \bigr|\le g(\omega)\qquad
\mbox{for all } (\omega,t)\in\Omega\times[0,t_*).
\end{equation}
Then (\ref{eqF=30}) holds.
\end{lemma}

Lemma~\ref{lemunderint,simpl} is
apparently rather common;
cf., for example, Theorem\break (2.27)(b) in \cite{folland}.
Lemma~\ref{lemunderint} will be used in the proof of Proposition~\ref
{propcard=1}. More generally, this lemma should be useful in certain
situations when the condition~(\ref{eq<g,simpl}) of the boundedness of
$(\partial_2f) (\omega,t)$ in $t$ for each $\omega$ is violated.
More specifically, in such situations (i) $(\partial_2f) (\omega,t)$
could have blow-up singularities and hence be unbounded in $t$ for each
$\omega$ in a somewhat ``small'' exceptional set $\Omega\setminus
\Omega_\varepsilon$, and yet (ii) the integration of $|(\partial_2f)
(\omega,t)|$ with respect to $|\mu(\dd\omega)|$ would smooth out
the singularities, resulting in a small value of the integral over the
``small'' set $\Omega\setminus\Omega_\varepsilon$---as is assumed
in (\ref{eqsupto0}).
Even though such situations seem rather natural and their treatment is
rather straightforward,
I have been unable to find in the literature a statement similar enough
to Lemma~\ref{lemunderint}.
So, for the readers' convenience, a
proof of Lemma~\ref{lemunderint} is provided below.

\begin{pf*}{Proof of Lemma~\ref{lemunderint}}
Take any $t\in[0,t_*)$,
$\varepsilon\in(0,\infty)$, and $h\in(0,h_{t,\varepsilon})$. Then
\begin{eqnarray}
\label{eq=I1+I2} \qquad\frac{F(t+h)-F(t)}h&=&\int_{\Omega}\mu(\dd\omega)
\frac{f(\omega,t+h)-f(\omega,t)}h\nonumber
\nonumber\\[-8pt]\\[-8pt]
&=&\int_{\Omega}\mu(\dd\omega) \int_0^1
\dd s (\partial_2f) (\omega,t+sh)
= I_{1,\varepsilon,h}(t)+I_{2,\varepsilon,h}(t),\nonumber
\end{eqnarray}
where
\[
I_{1,\varepsilon,h}(t):= \int_{\Omega_{t,\varepsilon}}\mu(\dd\omega)
\int
_0^1\dd s (\partial_2f) (
\omega,t+sh)
\]
and
\[
I_{2,\varepsilon,h}(t):= \int_{\Omega\setminus\Omega_{t,\varepsilon}}\mu
(\dd\omega) \int
_0^1\dd s (\partial_2f) (
\omega,t+sh).
\]
In view of the right continuity of $(\partial_2f) (\omega,t)$ in $t$
and the condition (\ref{eq<g}), by the dominated convergence theorem,
%
\begin{equation}
\label{eqI1to} I_{1,\varepsilon,h}(t)\underset{h\downarrow0}
{\longrightarrow}
I_{1,\varepsilon,0}(t)=\int_{\Omega_{t,\varepsilon}}\mu(\dd\omega) (
\partial_2f) (\omega,t).
\end{equation}
It also follows that the integral $I_{1,\varepsilon,0}(t)$ exists in
the Lebesgue sense (and is finite).
Next,
%
\begin{equation}
\label{eqsupI2to0} \sup_{h\in[0,h_{t,\varepsilon})} \bigl|I_{2,\varepsilon,h}(t)\bigr|\le\mathop{
\sup_{h\in[0,h_{t,\varepsilon})}}_{s\in(0,1)} \int_{\Omega\setminus
\Omega_{t,\varepsilon}} \bigl|
\mu(\dd\omega) (\partial_2f) (\omega,t+sh) \bigr| \underset{\varepsilon
\downarrow0} {\longrightarrow}0
\end{equation}
by (\ref{eqsupto0}). It also follows from (\ref{eqsupto0}) that the integral
\[
I_{2,\varepsilon,0}(t)=\int_{\Omega\setminus\Omega_{t,\varepsilon
}}\mu(\dd\omega) (
\partial_2f) (\omega,t)
\]
exists in the Lebesgue sense (and is finite) provided that $\varepsilon
$ is small enough.
So, the integral $I(t)$, defined in (\ref{eqF=30}), exists in the
Lebesgue sense (and is finite), since
$
I(t)=I_{1,\varepsilon,0}(t)+I_{2,\varepsilon,0}(t)$. 
Moreover, (\ref{eqsupI2to0}) implies
$
I_{2,\varepsilon,0}(t)
\underset{\varepsilon\downarrow0}{\longrightarrow}0$. %
Hence, $I_{1,\varepsilon,0}(t)=I(t)-I_{2,\varepsilon,0}(t)\underset
{\varepsilon\downarrow0}{\longrightarrow} I(t)$.
Combining now (\ref{eq=I1+I2}), (\ref{eqI1to}) and (\ref
{eqsupI2to0}), one completes the proof of the lemma.
\end{pf*}

\subsection{A calculus of variations of moments of infinitely divisible distributions with respect to
variations of the L\'evy characteristics}\label{var}

For any finite nonnegative Borel measure $H$ on $\mathbb{R}$,
let $Y_H$ denote any r.v. such that
%
\begin{equation}
\label{eqcf} \E e^{itY_H}=
\exp\biggl\{-t^2 \int
_\mathbb{R}H(\dd u) (R_1\exp) (0;itu) \biggr\}
\end{equation}
for all $t\in\mathbb{R}$,
where, for any $m\in\{0,1,\ldots\}$, any $(m+1)$-times continuously
differentiable function $g\dvtx \mathbb{R}\to\mathbb{R}$, and any
real $x$ and $u$,
%
\begin{eqnarray}
\qquad (R_m g) (x;u)&:=& \cases{ \displaystyle\frac{1}{u^{m+1}}
\Biggl(g(x+u)-\sum_{j=0}^m
\frac{u^j}{j!} g^{(j)}(x) \Biggr), &\quad if $u\ne0$,
\vspace*{3pt}\cr
\displaystyle\frac{1}{(m+1)!} g^{(m+1)}(x), &\quad if $u=0$}
\label{eqR==}
\\
&=&\frac{1}{m!} \int_0^1\dd s
(1-s)^m g^{(m+1)}(x+su). \label{eqR=}
\end{eqnarray}
This definition of $Y_H$ is valid, as the right-hand side expression in
(\ref{eqcf}) does define a characteristic function (c.f.) of (an
infinitely divisible) probability distribution, which is the weak limit
of a sequence of centered compound Poisson distributions.
Let
\[
\binom zj:=\frac{z(z-1)\cdots(z-j+1)}{j!}\qquad\mbox{for all }z\in
\mathbb{C}\mbox{ and all } j\in\{0,1,\ldots\}.
\]

%
\begin{lemma}\label{lemdiff}
Take any $q\in(2,\infty)$, $t_0\in(0,\infty)$ and $\sigma\in
(0,\infty)$.
Let $H$ be a nonnegative Borel measure on $\mathbb{R}$, and let
$\Delta$ be a real-valued Borel measure on $\mathbb{R}$ such that the measure
%
\begin{equation}
\label{eqHt=} H_t:=H+t\Delta
\end{equation}
is nonnegative for all $t\in[0,t_0]$.
Let $Y$ be a r.v. independent of $Y_{H_t}$, where $Y_{H_t}$ is defined
according to (\ref{eqcf}).
Suppose also that
%
\begin{equation}
\label{eqexpmom} \int_\mathbb{R} \bigl(\P(Y\in\dd u)+H(\dd u)+\bigl|
\Delta(\dd u)\bigr| \bigr) e^{\sigma|u|}<\infty.
\end{equation}
Then for all $t\in[0,t_0)$
%
\begin{eqnarray}\label{eqder1}
&& \biggl(\pd{} {t} \biggr)^+\E|Y+Y_{H_t}|^q
\nonumber\\[-8pt]\\[-8pt]
&&\qquad =2!
\binom q2 \int_\mathbb{R}\Delta(\dd u)\int_0^1
\dd s (1-s) \E|su+Y+Y_{H_t}|^{q-2},\nonumber
\end{eqnarray}
where
$ (\pd{}{t} )^+$ denotes the right-hand side partial
derivative in $t$.

Moreover, if $q>4$, then
for all $t\in[0,t_0)$
%
\begin{eqnarray}\label{eqder2}
&& \biggl(\pd{{}^2} {t^2} \biggr)^+ \E|Y+Y_{H_t}|^q\nonumber
\\
&&\qquad =4! \binom q4
\int_{\mathbb{R}^2}\Delta(\dd
u_1)\Delta(\dd u_2)
\\
&&\quad\qquad{}\times  \int_{(0,1)^2} \dd
s_1\,\dd s_2 (1-s_1) (1-s_2) \E
|s_1u_1+s_2u_2+Y+Y_{H_t}|^{q-4},\nonumber
\end{eqnarray}
where $ (\pd{{}^2}{t^2} )^+:= (\pd{}{t} )^+ (\pd
{}{t} )^+$ denotes the second right-hand side partial derivative in $t$.

Identities (\ref{eqder1}) and (\ref{eqder2}) hold if the four
instances therein of the absolute-value function $x\mapsto|x|$ are
replaced by the four instances of the positive-part function $x\mapsto
x_+:=0\vee x$ or the four instances of the negative-part function
$x\mapsto x_-:=(-x)_+$.
\end{lemma}

\begin{pf}
By Theorem~1 in \cite{positive},
%
\begin{equation}
\label{eqstart} \E(Y+Y_{H_t})_+^q=\kappa_q\int
_{\Re z=\sigma}\frac{\dd z}{z^{q+1}} \E e^{zY} \E\exp
\{zY_{H_t}\},
\end{equation}
where $x_+^q:=(x_+)^q$ for all $x\in\mathbb{R}$ and
%
\begin{equation}
\label{eqkaq} \kappa_q:=\frac{\Gamma(q+1)}{2\pi i};
\end{equation}
here and below in this proof, by default, $t\in[0,t_0)$.
By (\ref{eqcf}) and analytic continuation, for all $z\in\mathbb{C}$
with $\Re z=\sigma$
%
\begin{equation}
\label{eqcf,z} \E\exp\{zY_{H_t}\}=\exp\biggl\{z^2 \int
_\mathbb{R}H_t(\dd u) (R_1\exp) (0;zu)
\biggr\},
\end{equation}
whence, by (\ref{eqHt=}),
%
\begin{equation}
\label{eqD+E,2} \biggl(\pd{} {t} \biggr)^+ \E\exp\{zY_{H_t}\} =\E\exp
\{zY_{H_t}\} z^2 \int_\mathbb{R}\Delta(\dd
u) (R_1\exp) (0;zu).
\end{equation}
In view of (\ref{eqR==}), $\Re[z^2(R_1\exp)(0;zu) ]\le
\sigma
^2(R_1\exp)(0;\sigma u)\le\sigma^2e^{\sigma|u|}/2$ and $|(R_1\exp
)(0;zu)|\le e^{\sigma|u|}/2$ for all $u\in\mathbb{R}$ and all $z\in
\mathbb{C}$ with $\Re z=\sigma>0$.
It follows by (\ref{eqcf,z}) and (\ref{eqexpmom}) that, again for all
$z\in\mathbb{C}$ with $\Re z=\sigma>0$,
\[
\sup_{t\in[0,t_0]}\bigl|\E\exp\{zY_{H_t}\}\bigr| \le\exp\biggl\{
\frac{\sigma^2}2 \int_\mathbb{R} \bigl(H(\dd
u)+t_0\bigl|\Delta(\dd u)\bigr| \bigr) e^{\sigma|u|} \biggr\}<\infty
\]
and
\[
\biggl|\int_\mathbb{R}\Delta(\dd u) (R_1\exp) (0;zu) \biggr|
\le\int_\mathbb{R}\bigl|\Delta(\dd u)\bigr| e^{\sigma|u|}/2<\infty.
\]
Also, $\int_{\Re z=\sigma} |z^2\frac{\dd z}{z^{q+1}}
|<\infty
$, since $q>2$ and $\sigma\in(0,\infty)$.
So, by Lemma~\ref{lemunderint,simpl}, 
%
\begin{equation}
\label{eqD+E,1} \biggl(\pd{} {t} \biggr)^+\E(Y+Y_{H_t})_+^q =
\kappa_q\int_{\Re z=\sigma}\frac{\dd z}{z^{q+1}} \E
e^{zY} \biggl(\pd{} {t} \biggr)^+\E\exp\{zY_{H_t}\}.
\end{equation}
Further, by (\ref{eqR=}),
$(R_1\exp)(0;zu)=\int_0^1\dd s (1-s) e^{szu}$.
Hence, by (\ref{eqD+E,1}), (\ref{eqD+E,2}), the Fubini theorem, (\ref
{eqkaq}), and (\ref{eqstart}),
\begin{eqnarray*}
&& \biggl(\pd{} {t} \biggr)^+ \E(Y+Y_{H_t})_+^q
\\
&&\qquad = \kappa_q\int_{\Re z=\sigma}\frac{\dd z}{z^{q+1}} \E
e^{zY} \E e^{zY_{H_t}} z^2 \int_\mathbb{R}
\Delta(\dd u) \int_0^1\dd s (1-s)
e^{szu}
\\
&&\qquad =\frac{\kappa_q}{\kappa_{q-2}} \int_\mathbb{R}\Delta(\dd u) \int
_0^1\dd s (1-s) \kappa_{q-2}\int
_{\Re z=\sigma}\frac{\dd
z}{z^{q-1}} \E e^{z(su+Y)} \E
e^{zY_{H_t}}
\\
&&\qquad = 2!\binom q2 \int_\mathbb{R}\Delta(\dd u) \int
_0^1\dd s (1-s) \E(su+Y+Y_{H_t})_+^{q-2}.
\end{eqnarray*}
This proves (\ref{eqder1}) for the function $x\mapsto x_+$ in place of
the function $x\mapsto|x|$.
Now~(\ref{eqder2}), again for the function $x\mapsto x_+$, follows by
Lemma~\ref{lemunderint,simpl}. 

The case of the function $x\mapsto x_-$ can be considered quite
similarly. Alternatively, this case can be simply reduced to the case
of the function $x\mapsto x_+$ by observing that, with $H^-_t(\dd
u):=H_t(-\dd u)$, one has $-Y_{H_t}\stackrel{\mathrm{D}}=Y_{H^-_t}$,
where $\stackrel{\mathrm{D}}=$ denotes the equality in distribution.

Finally, the case of the function $x\mapsto|x|$ follows immediately
from the considered two cases by the obvious identity
$|x|^r=x_+^r+x_-^r$ for all $r\in(0,\infty)$ and $x\in\mathbb{R}$.
\end{pf}

Results similar to Lemma~\ref{lemdiff}, but for general moment
functions $f$ in place of the power-like moment functions $|\cdot|^q$,
$\cdot_+^q$ and $\cdot_-^q$ in Lemma~\ref{lemdiff}, were obtained by
lengthier direct probabilistic arguments in earlier versions of this
paper \cite{rosenthalbestv2}.
It is possible to obtain such more general results by the
Fourier--Laplace method as well, by decomposing $f$ into\vspace*{1pt} harmonics, the
way this was done in \cite{positive} for the function $\cdot_+^p$.
However, this possibility will not be pursued here.

\subsection{Main propositions in the proof of Theorem~\texorpdfstring{\protect\ref{th}}{1.3}}\label{mainprops}

Let $\mathscr{H}$ denote the set of all nonnegative Borel measures on
$\mathbb{R}$.
Take any real numbers $p>3$, $A>0$, $B>0$, and $M>0$, and introduce the
following subsets of the set $\mathscr{H}$:
%
\begin{eqnarray}
\mathscr{H}_{p;A,B}&:=&\biggl\{H\in\mathscr{H}\dvtx \int H(\dd x)=B,\int
|x|^{p-2}H(\dd x)=A\biggr\}, \label{eqHAB}
\\
\mathscr{H}_{p;\le A,\le B}&:=&\biggl\{H\in\mathscr{H}\dvtx \int H(\dd
x)\le B,
\int|x|^{p-2}H(\dd x)\le A\biggr\}, \label{eqH<A<B}
\\
\mathscr{H}_{p;A,B;M}&:=&\bigl\{H\in\mathscr{H}_{p;A,B}\dvtx \supp H
\subseteq[-M,M]\bigr\}, \label{eqHAB;M}
\\
\mathscr{H}_{p;\le A,\le B;M}&:=&\bigl\{H\in\mathscr{H}_{p;\le A,\le
B}\dvtx \supp
H\subseteq[-M,M]\bigr\}, \label{eqH<A<B;M}
\end{eqnarray}
where $\supp H$ stands for the support set of the measure $H$; we also
write $\int$ for $\int_\mathbb{R}$.
Note that the set $\mathscr{H}_{p;\le A,\le B}$ obviously contains the other
three of the above four sets.

%
\begin{remark}\label{remnonempty}
Given any positive real $A$, $B$, and $M$,
for the condition $\mathscr{H}_{p;A,B;M}\ne\varnothing$ to hold it is clearly
necessary that
%
\begin{equation}
\label{eqA<} A\le BM^{p-2}
\end{equation}
or, equivalently, $B\ge A/M^{p-2}$ or, equivalently,
%
\begin{equation}
\label{eqM>} M\ge c, 
\end{equation}
where $c=c_p(A,B)$ as in (\ref{eqla,c}).

Therefore, in the statements concerning
$\mathscr{H}_{p;A,B;M}$, let us assume by default that this
restriction on $A$,
$B$, and $M$ holds.
\end{remark}

\begin{center}
\framebox{In Propositions \ref{propattain}--\ref{propcard=1} below,
let $X$ be any bounded zero-mean r.v.}
\end{center}

Let then
%
\begin{eqnarray}
\mathscr{S}_{p,q;A,B;X;M}&:=&\sup\bigl\{\E|X+Y_H|^q
\dvtx  H\in\mathscr{H} _{p;A,B;M}\bigr\}, \label{eqSAB;M}
\\
\mathscr{S}_{p,q;\le A,\le B,X;M}&:=&\sup\bigl\{\E|X+Y_H|^q
\dvtx  H\in\mathscr{H} _{p;\le A,\le B;M}\bigr\}, \label{eqS<A<B;M}
\end{eqnarray}
where $Y_H$ and $q$ are as in (\ref{eqcf}) and (\ref{eqq}), respectively.

%
\begin{proposition}\label{propattain}
The supremum $\mathscr{S}_{p,q;\le A,\le B,X;M}$ is finite and attained.
If $A\le BM^{p-2}$ (recall Remark~\ref{remnonempty}), then the
supremum $\mathscr{S}_{p,q;A,B;X;M}$ is finite and attained as well.
\end{proposition}

%
\begin{proposition}\label{propincr}
Suppose that $p\ge q>4$.
Then the supremum\break
$\mathscr{S}_{p,q;A,B;X;M}$ is (strictly) increasing in $B\in
[A/M^{p-2},\infty)$ for each $A\in(0,\infty)$ and in $A\in
(0,BM^{p-2}]$ for each $B\in(0,\infty)$;
in particular, it follows that
%
\begin{equation}
\label{eqle=eq} \mathscr{S}_{p,q;A,B;X;M}=\mathscr{S}_{p,q;\le A,\le B,X;M}
\end{equation}
for any positive real $A$, $B$, and $M$ such that $A\le BM^{p-2}$.
\end{proposition}

Introduce the set
%
\begin{eqnarray}\label{eqH*}
&& \mathscr{H}_{*,p,q;A,B;X;M}
\nonumber\\[-8pt]\\[-8pt]
&&\qquad:= \bigl\{H\in
\mathscr{H}_{p;A,B;M}\dvtx \E|X+Y_H|^q=\mathscr
{S}_{p,q;A,B;X;M}=\mathscr{S}_{p,q;\le A,\le B,X;M} \bigr\}\hspace*{-18pt}\nonumber
\end{eqnarray}
of the maximizers of $\E|X+Y_H|^q$ over all $H\in\mathscr
{H}_{p;A,B;M}$ or,
equivalently, over all $H\in\mathscr{H}_{p;\le A,\le B;M}$.
According to Propositions~\ref{propattain} and \ref{propincr},
\[
\mathscr{H}_{*,p,q;A,B;X;M}\ne\varnothing.
\]

%
\begin{proposition}\label{prop<2}
Suppose that $p\ge q>5$. Take any $H\in\mathscr{H}_{*,p,q;A,B;X;M}$.
Then
%
\begin{equation}
\label{eqle1} \card\bigl((0,\infty)\cap\supp H \bigr)\le1\quad\mbox
{and}\quad
\card\bigl((-\infty,0)\cap\supp H \bigr)\le1,
\end{equation}
where $\card$ denotes the cardinality of the set.
\end{proposition}

%
\begin{proposition}\label{propnogaus}
Suppose that $p\ge q>4$. Take any $H\in\mathscr{H}_{*,p,q;A,B;X;M}$.
Then
$H(\{0\})=0$.
\end{proposition}

%
\begin{proposition}\label{propcard=1}
Suppose that $p\ge q>5$. Take any $H\in\mathscr{H}_{*,p,q;A,B;X;M}$.
Suppose also that the set $\supp H$ is contained in the \emph{open}
interval $(-M,M)$.
Then $\card\supp H=1$.
\end{proposition}

%
\begin{proposition}\label{proplim}
Suppose that $p\ge q>2$.
Let the quadruple $(c_1,c_2,\break w_1,w_2)\in\mathbb{R}^2\times[0,\infty
)^2$ vary so\vspace*{1pt} that
$w_1+w_2=B$,
$c_1\to b$, $|c_2|\to\infty$, and $|c_2|^{q-2}w_2\to a$, for some
$a\in[0,A]$ and $b\in[-c,c]$, where $c$ is as in (\ref{eqla,c}).
Then, for $H:=H_{c_1,c_2,w_1,w_2}:=w_1\delta_{c_1}+w_2\delta_{c_2}$,
\[
\E|X+Y_H|^q\longrightarrow a
+
\E|X+Y_{B\delta_b}|^q.
\]
\end{proposition}

\begin{pf*}{Proof of Proposition~\ref{propattain}}
Let us only show that the supremum\break $\mathscr{S}_{p,q;A,B;X;M}$ is
finite and attained; that $\mathscr{S}_{p,q;\le A,\le B,X;M}$ is so is
shown similarly and even a bit more easily.
Let $(H_m)$ be a sequence in $\mathscr{H}_{p;A,B;M}$ such that $\E
|X+Y_{H_m}|^q\to\mathscr{S}_{p,q;A,B;X;M}$.
Because the interval $[-M,M]$ is compact and the functions $1$ and
$|\cdot|^{p-2}$ are continuous and bounded on $[-M,M]$, without loss
of generality (w.l.o.g.) the sequence $(H_m)$ converges weakly to some
$H\in\mathscr{H}_{p;A,B;M}$. So, by (\ref{eqcf}) and (\ref{eqR=}),
$Y_{H_m}\stackrel{\mathrm{D}}\longrightarrow Y_H$, since $(R_1\exp
)(0;itu)$ is continuous and bounded
in $u\in[-M,M]$.
Moreover, by the analytic extension of (\ref{eqcf}), for any $\tH\in
\mathscr{H}_{p;A,B;M}$
%
\begin{eqnarray}\label{eqcosh}
\E\cosh(kY_\tH) &=& \frac12 \exp\biggl
\{k^2 \int\tH(\dd u) (R_1\exp) (0;ku) \biggr\}\nonumber
\\
&&{}  +\frac12 \exp\biggl\{k^2 \int\tH(\dd u) (R_1\exp) (0;-ku) \biggr\}
\\
&\le& \exp\bigl\{k^2B (R_1\exp) \bigl(0;|k|M\bigr)\bigr\}<\infty\nonumber
\end{eqnarray}
for all real $k$---because, by (\ref{eqR=}), $(R_1\exp)(0;u)$ is increasing in $u\in
\mathbb{R}$.
Also, $|X+Y_{H_m}|^q\le2^{q-1}(|X|^q+|Y_{H_m}|^q)$. 
So, by \cite{billingsley}, Theorem~5.4,
%
\begin{equation}
\label{eqattain} \mathscr{S}_{p,q;A,B;X;M}=\lim_m
\E|X+Y_{H_m}|^q=\E|X+Y_H|^q<\infty.
\end{equation}\upqed
\end{pf*}

\begin{pf*}{Proof of Proposition~\ref{propincr}}
Let us show that $\mathscr{S}_{p,q;A,B;X;M}$ is increasing in $A$ and
in $B$;
then (\ref{eqle=eq}) follows immediately.

In accordance with 
Proposition~\ref{propattain}, take any $H\in\mathscr{H}_{p;A,B;M}$
such that
$\E|X+Y_H|^q=\mathscr{S}_{p,q;A,B;X;M}$.
Then
$H_t:=H+t\delta_0\in\mathscr{H}_{p;A,B+t;M}$
for all real $t\ge0$, where, as before, $\delta_u$ denotes the Dirac
probability measure at $u$.
So, by Lemma~\ref{lemdiff}, the right derivative of $\E|X+Y_{H_t}|^q$
in $t$ at $t=0$ is
${q\choose 2} \E|X+Y_H|^{q-2}>0$; the\vspace*{1pt} last inequality is strict because
the measure $H$ is in $\mathscr{H}_{p;A,B;M}$ and hence nonzero, which
in turn
implies that the r.v. $Y_H$ is nondegenerate.
Therefore, for the lower right derivative of $\mathscr
{S}_{p,q;A,B;X;M}$ in $B$ one has
\begin{eqnarray*}
\liminf_{t\downarrow0}\frac{\mathscr{S}_{p,q;A,B+t,X;M}-\mathscr
{S}_{p,q;A,B;X;M}}t &\ge&\liminf
_{t\downarrow0}\frac{\E|X+Y_{H_t}|^q-\E|X+Y_H|^q}t
\\
&=&\lim_{t\downarrow0}\frac{\E|X+Y_{H_t}|^q-\E|X+Y_H|^q}t
\\
&=&\binom q2 \E
|X+Y_H|^{q-2}>0.
\end{eqnarray*}
Next, note that $\mathscr{S}_{p,q;A,B;X;M}$ is left-upper semi-continuous
in $B\in(A/M^{p-2},\break \infty)$; that is,
\[
\limsup_{\tilde{B}\uparrow B}\mathscr{S}_{p,q;A,\tilde{B},X;M}\le
\mathscr{S}_{p,q;A,B;X;M}.
\]
Indeed, take any sequence $(B_m)$ such that $B_m\uparrow B$ and
\[
\lim_{m\to\infty}\mathscr{S}_{p,q;A,B_m,X;M}>\mathscr{S}_{p,q;A,B;X;M}.
\]
By Proposition~\ref{propattain}, for each large enough $m$ there is
some measure $H_m\in\mathscr{H}(p,A,B_m;M)$ such that
$\E|X+Y_{H_m}|^q=\mathscr{S}_{p,q;A,B_m,X;M}$.
Passing to a subsequence of the sequence $(B_m)$,
w.l.o.g. one may assume that $H_m$ converges weakly on the compact set
$[-M,M]$ to some measure $H_*$.
Since the functions~$1$, $|\cdot|^{p-2}$, and
$|\cdot|^q$ are continuous, it follows that $H_*\in\mathscr
{H}(p,A,B;M)$ and
$\mathscr{S}_{p,q;A,B_m,X;M}=\E|X+Y_{H_m}|^q\longrightarrow\E|X+Y_{H_*}|^q
\le\mathscr{S}_{p,q;A,B;X;M}$ as $m\to\infty$, which contradicts
the assumption on the sequence $(B_m)$.
This completes the proof that $\mathscr{S}_{p,q;A,B;X;M}$ is
increasing in $B$.

To show that $\mathscr{S}_{p,q;A,B;X;M}$ is increasing in $A$, take any
$A\in(0,BM^{p-2})$; cf.~(\ref{eqA<}).
Then
%
\begin{equation}
\label{eqH>0} H \bigl((-M,M) \bigr)>0,
\end{equation}
because otherwise $\supp H\subseteq\{-M,M\}$ and hence $A=BM^{p-2}$.
So, there exists some $b\in(-M,M)\cap\supp H$.
For $\delta\in(0,\infty)$ and $t\in[0,\infty)$, let now
%
\begin{equation}
\label{eqHt,de} H_t:=H_{\delta,t}:=H+t\Delta,
\end{equation}
where $\Delta=\Delta_\delta$ is the real-valued Borel measure on
$\mathbb{R}$ defined by the condition that
%
\begin{eqnarray}
\label{eqDede}
\int_\mathbb{R}f(u)\Delta(\dd u)
&=& \frac{M+b}{2M} f(M)+\frac{M-b}{2M} f(-M)
\nonumber\\[-8pt]\\[-8pt]
&&{} -\frac{1}{H([b-\delta,b+\delta])} \int
_{[b-\delta,b+\delta]}f(u) H(\dd u)\nonumber
\end{eqnarray}
for all locally bounded (say) Borel functions $f\dvtx \mathbb{R}\to
\mathbb{R}$;
note that $H([b-\delta,b+\delta])>0$, by the condition $b\in\supp H$.
Also, then the measure $H_t$
is nonnegative for all $t\in[0,t_0]$, where $t_0:=H([b-\delta,b+\delta])>0$.
So, letting
%
\begin{equation}
\label{eqhq} h(x):=2! \binom q2 \E|x+X+Y_H|^{q-2}
\end{equation}
for all $x\in\mathbb{R}$,
by Lemma~\ref{lemdiff} one has
%
\begin{eqnarray}\label{eqder}
&& \lim_{\delta\downarrow0}\lim_{t\downarrow0}\frac{\E
|X+Y_{H_t}|^q-\E|X+Y_H|^q}t
\nonumber\\[-8pt]\\[-8pt]
&&\qquad =\int_0^1\dd s
(1-s) \biggl(\frac{M+b}{2M}h(Ms)+\frac
{M-b}{2M}h(-Ms)-h(bs) \biggr)>0,\nonumber
\end{eqnarray}
because $q>4$, and the r.v. $Y_H$ is nondegenerate, whence the function
$h$ is strictly convex.
Thus,
eventually
%
\begin{equation}
\label{eq<} \mathscr{S}_{p,q;A,B;X;M}=\E|X+Y_H|^q<
\E|X+Y_{H_t}|^q.
\end{equation}
In this context, we say that an assertion $\mathscr{A}=\mathscr
{A}_{\delta,t}$ holds ``eventually'' if
$\exists\delta_*\in(0,\infty)$
$\forall\delta\in(0,\delta_*)$ $\exists t_\delta\in(0,t_0)$
$\forall t\in(0,t_\delta)$ $\mathscr{A}_{\delta,t}$ holds; recall
here that, in view of~(\ref{eqHt,de}) and~(\ref{eqDede}), $H_t$
depends not only on $t$ but also on $\delta$.

On the other hand, for all $t\in(0,t_0)$ one has $\int_\mathbb
{R}H_t(\dd x)=\int_\mathbb{R}H(\dd x)+t\int_\mathbb{R}\Delta(\dd
x)=B+t\int_\mathbb{R}\Delta(\dd x)=B$ and $\int_\mathbb
{R}|x|^{p-2}H_t(\dd x)=A+ta$ and hence $H_t\in\mathscr
{H}(p,A+ta,B;M)$, where
$a:=\int_\mathbb{R}|x|^{p-2}\Delta(\dd x)\ge
(M^{p-2}-(|b|+\delta)^{p-2} )>0$ for all small enough $\delta>0$.
So, by (\ref{eqSAB;M}), eventually $\mathscr{S}_{p,q;A+ta,B,X;M}\ge
\E|X+Y_{H_t}|^q$, whence, by (\ref{eq<}), $\mathscr{S}_{p,q;\cdot,B,X;M}$ is increasing in a right neighborhood of the previously chosen
value of $A\in(0,BM^{p-2})$.

Since $A$ was chosen arbitrarily in the interval $(0,BM^{p-2})$, to
complete the proof of Proposition~\ref{propincr},
it remains to note that $\mathscr{S}_{p,q;A,B;X;M}$ is left-upper
semi-continuous
in $A\in(0,BM^{p-2}]$; this semi-continuity property is established
quite similarly to the left-upper semi-continuity in $B$, proved earlier.
\end{pf*}

\begin{pf*}{Proof of Proposition~\ref{prop<2}}
To obtain a contradiction, suppose that there exist $b$ and $b_1$ such
that $0<b<b_1<\infty$ and $\{b,b_1\}\subseteq\supp H$.
In view of possible rescaling [i.e., replacing $X$, $A$, $B$, $M$
and $H(\dd x)$ by $X/b_1$, $A/b_1^p$, $B/b_1^2$, $M/b_1$,
and $H(b_1 \,\dd y)/b_1^2$, resp.],
w.l.o.g. assume that $b_1=1$, so that
\[
0<b<1.
\]
By (\ref{eqH*}),
%
\begin{equation}
\label{eqatta} \mathscr{S}_{p,q;\le A,\le B,X;M} =\E|X+Y_H|^q.
\end{equation}

Introduce now
%
\begin{equation}
\label{eqk} k:=\frac{b^2(1-b^{p-3})}{p-3},
\end{equation}
take any
%
\begin{equation}
\label{eqa} a\in(0,1/k)
\end{equation}
and then also introduce
%
\begin{eqnarray}\label{eqvp,ta,tb}
\varepsilon &:=& a \bigl(b-b^{p-1}-(p-2)k \bigr),\nonumber
\\
\tilde{a}&:=&1+a \bigl(b^{p-1}+(p-1)k \bigr), \quad\mbox{and}
\\
\tilde{b}&:=&1-ka.\nonumber
\end{eqnarray}
Note that the conditions (\ref{eqvp,ta,tb}) and (\ref{eqa}) imply
${\tilde{b}}\in(0,1)$.
Observe also that $\varepsilon=abr(b)/(p-3)$, where
$r(b):=p-3-(p-2)b+b^{p-2}$, and $r(1)=0$ and $r'(b)=-(p-2)(1-b^{p-3})<0$
for $b\in(0,1)$, so that $\varepsilon>0$.

Define the real-valued measure $\Delta=\Delta_{a,\delta}$ by the condition
%
\begin{eqnarray}\label{eqDe}
&& \int_\mathbb{R}f(u)\Delta(\dd u)\nonumber
\\
&&\qquad =\varepsilon
f(0)+{\tilde{a}} {\tilde{b}}f({\tilde{b}}) -\frac{ab}{H([b-\delta,b+\delta])} \int
_{[b-\delta,b+\delta
]}f(u)H(\dd u)
\\
&&\quad\qquad{}  -\frac{1}{H([1-\delta,1+\delta])} \int_{[1-\delta,1+\delta
]}f(u)H(\dd u)\nonumber
\end{eqnarray}
for all locally bounded (say) Borel functions $f\dvtx \mathbb{R}\to
\mathbb{R}$, where
$\delta$ is any real number in the interval $(0,\frac{1-b}2)$, so
that the denominators $H([b-\delta,b+\delta])$ and $H([1-\delta,1+\delta
])$ are strictly positive, and the intervals $[b-\delta,b+\delta]$ and
$[1-\delta,1+\delta]$ are disjoint, in view of the
assumptions $\{b,b_1\}\subseteq\supp H$ and $b_1=1$.
For
$t\in[0,\infty)$, let now
%
\begin{equation}
\label{eqHt,a,de} H_t:=H_{a,\delta,t}:=H+t\Delta.
\end{equation}
This measure is nonnegative for all $t\in[0,t_0]$, where
\[
t_0:=\min\biggl(\frac1{ab} H\bigl([b-\delta,b+\delta]\bigr),H
\bigl([1-\delta,1+\delta]\bigr) \biggr)>0.
\]
By Lemma~\ref{lemdiff},
\begin{eqnarray*}
&& \lim_{\delta\downarrow0}\lim_{t\downarrow0}\frac{\E
|X+Y_{H_{a,\delta,t}}|^q-\E|X+Y_H|^q}t
\\
&&\qquad =\int_0^1\dd s (1-s) \bigl[
\varepsilon h(0)+{\tilde{a}} {\tilde{b}}h(s{\tilde{b}})-ab h(sb)-h(s)
\bigr],
\end{eqnarray*}
where the function $h$ is still defined by (\ref{eqhq}).

Letting further $a\downarrow0$ and using Lemma~\ref{lemunderint,simpl}, one obtains
%
\begin{eqnarray}\label{eqL=}
\mathscr{L}&:=&\frac{1}{b^2} \lim_{a\downarrow0}
\frac{1}a \lim_{\delta
\downarrow0} \lim_{t\downarrow0}
\frac{\E|X+Y_{H_{a,\delta,t}}|^q-\E|X+Y_H|^q}t
\nonumber\\[-8pt]\\[-8pt]
& =&\int_0^1\dd s (1-s)F(b,s),\nonumber
\end{eqnarray}
where, in view of 
(\ref{eqvp,ta,tb}), (\ref{eqk}),
(\ref{eqR==}) and (\ref{eqR=}),
%
\begin{eqnarray}
\label{eqF=77} F(b,s)&:=&\frac{1}{b^2}\frac{\dd}{\dd a} \bigl[\varepsilon
h(0)+{\tilde{a}} {\tilde{b}}h(s{\tilde{b}})-ab h(sb)-h(s) \bigr]
\Big|_{a=0}\nonumber
\\
&=&\frac{1}{b^2} \bigl\{
\bigl(b^{p-1}+(p-2)k \bigr)
\bigl(h(s)-h(0) \bigr)
-b \bigl(h(sb)-h(0) \bigr)-ks
h'(s) \bigr\}\nonumber
\\
&=&h(s)-h(0)-\frac{h(sb)-h(0)}b +\frac{1-b^{p-3}}{p-3} \bigl[h(s)-h(0)-sh'(s)
\bigr]
\\
&=&\int_0^1\dd\alpha\, s \bigl[h'(
\alpha s)-h'(\alpha sb) \bigr] -\int_b^1
\dd u\, u^{p-4} s^2(R_1h) (s;-s)\nonumber 
\\
&=&\int_0^1\dd\alpha\, s \alpha s\int_b^1\dd u\, h''(u
\alpha s) -\int_b^1\dd u\, u^{p-4}
s^2 \int_0^1\dd\theta(1-\theta)
h''(s-\theta s)\nonumber
\\
&=&s^2\int_b^1\dd u \int
_0^1\dd\alpha\,\alpha\bigl[h''(u
\alpha s)-u^{p-4}h''(\alpha s) \bigr].\nonumber
\end{eqnarray}
%

By (\ref{eqhq}) and Lemma~\ref{lemunderint,simpl}, for $x\in\mathbb
{R}$ and $u\in(0,\infty)$
\[
\frac{u^{q-4}h''(x)}{q(q-1)(q-2)(q-3)}=\psi_{ux}(u),
\]
where $\psi_v(u):=\E|v+uW|^{q-4}$ for all $v\in\mathbb{R}$ and $W:=X+Y_H$.
Note that, in view of the condition $\E X=0$ and the definition (\ref{eqcf}),
$\E W=0$. Also, $\E W^2>0$,
because $\card\supp H\ge2>0$ and hence $H\ne0$ and thus the r.v.
$Y_H$ is nondegenerate.
Also, clearly $h''\ge0$.
Therefore and because $p\ge q>5$, for all $u$, $\alpha$ and $s$ in $(0,1)$
%
\begin{eqnarray}\label{eq>0}
\frac{h''(u\alpha s)-u^{p-4}h''(\alpha
s)}{q(q-1)(q-2)(q-3)} &\ge&\frac{h''(u\alpha s)-u^{q-4}h''(\alpha
s)}{q(q-1)(q-2)(q-3)} 
\nonumber\\[-8pt]\\[-8pt]
&=& \psi_{u\alpha s}(1)-\psi_{u\alpha s}(u).\nonumber
\end{eqnarray}
%
Recalling that $q>5$, one sees that for each $v\in(0,\infty)$ the
function $\psi_v$ is convex, with $\psi'_v(0)=0$ and
$\psi''_v(0)=(q-4)(q-5)\E W^2 v^{q-6}>0$.
This implies that $\psi_v(u)$ is strictly increasing in $u\ge0$,
which shows that the expression  $\psi_{u\alpha s}(1)-\psi_{u\alpha
s}(u)$ in (\ref{eq>0}) is strictly positive.
Thus, by (\ref{eqL=}), (\ref{eqF=77}) and (\ref{eq>0}), $\mathscr{L}>0$.
Now (\ref{eqatta}) implies that
eventually
\[
\mathscr{S}_{p,q;\le A,\le B,X;M}=\E|X+Y_H|^q<
\E|X+Y_{H_{a,\delta,t}}|^q.
\]
In this context, we say that an assertion $\mathscr{A}=\mathscr
{A}_{a,\delta,t}$ holds ``eventually'' if
$\exists a_0\in(0,\infty)$ $\forall a\in(0,a_0)$ 
$\exists\delta^*_a\in(0,\infty)$
$\forall\delta\in(0,\delta^*_a)$ $\exists t_{a,\delta}\in(0,t_0)$
$\forall t\in(0,t_{a,\delta})$ $\mathscr{A}_{a,\delta,t}$ holds.

Thus, we obtain a contradiction with the definition of $\mathscr
{S}_{p,q;\le A,\le B,X;M}$ in (\ref{eqS<A<B;M}), because, as we shall
check in moment, $H_{a,\delta,t}\in\mathscr{H}_{p;\le A,\le B;M}$ eventually.
Indeed, by (\ref{eqDe}), (\ref{eqvp,ta,tb}), and (\ref{eqk}),
\begin{eqnarray*}
\int_\mathbb{R}\Delta(\dd x)&=&\varepsilon+{
\tilde{a}} {\tilde{b}}-ab-1
\\
&=&a \bigl(b-b^{p-1}-(p-2)k \bigr)
\\
&&{} +(1-ka) \bigl(1+a\bigl[b^{p-1}+(p-1)k
\bigr] \bigr)-ab-1
\\
&<&ab-a \bigl(b^{p-1}+(p-2)k \bigr)+1+\bigl[b^{p-1}+(p-2)k
\bigr]a-ab-1
\\
&=&0,
\end{eqnarray*}
so that
%
\begin{equation}
\label{eq<B} \int_\mathbb{R}H_{a,\delta,t}(\dd x)=\int
_\mathbb{R}H(\dd x)+t\int_\mathbb{R}\Delta(\dd
x)<\int_\mathbb{R}H(\dd x)=B.
\end{equation}
%
Similarly,
\begin{eqnarray*}
&&\lim_{\delta\downarrow0}  \int_\mathbb{R}|x|^{p-2}
\Delta(\dd x)={\tilde{a}} {\tilde{b}}^{p-1}-ab^{p-1}-1<0,
\end{eqnarray*}
where the inequality holds eventually, for all small enough $a>0$.
Indeed, in view of (\ref{eqvp,ta,tb}), this inequality can be
rewritten as
%
\begin{equation}
\label{eqf<0} f_\gamma(u):=\bigl[1+(\gamma+r)u\bigr](1-u)^r-(1+
\gamma u)<0,
\end{equation}
with $r:=p-1>0$, $u:=ka$, and $\gamma:=b^r/k\ge0$. Note that
eventually $u\in(0,1)$.
To verify inequality (\ref{eqf<0}) for such $u$, note that $f_\gamma
(u)$ decreases in $\gamma$, so that w.l.o.g. $\gamma=0$. The
inequality $f_0(u)<0$ is equivalent to $\ln(1+ru)+r\ln(1-u)<0$, which
is easy to check for $u\in(0,1)$ by differentiation.
It follows that [cf. (\ref{eq<B})] $\int_\mathbb
{R}|x|^{p-2}H_{a,\delta,t}(\dd x)<\int_\mathbb{R}|x|^{p-2}H(\dd x)=A$
eventually.

Also, the conditions $H\in\mathscr{H}_{p;\le A,\le B;M}$, 
$\{b,b_1\}\subseteq\supp H$ and $b_1=1$ imply $\supp H\subseteq
[-M,M]$ and hence $M\ge1$.
So,
$\supp H_{a,\delta,t}\subseteq\supp H\cup\{0,{\tilde{b}}\}\subseteq
[-M,M]$ eventually, in view of (\ref{eqvp,ta,tb}).

By (\ref{eqH<A<B;M}), we conclude that indeed $H_{a,\delta,t}\in
\mathscr{H}
_{p;\le A,\le B;M}$ eventually.
Thus, indeed the assumption that there exist $b$ and $b_1$ such that
$0<b<b_1<\infty$ and $\{b,b_1\}\subseteq\supp H$ leads to a
contradiction, which proves the first inequality in (\ref{eqle1}).
The second inequality there can be proved quite similarly or,
alternatively, quickly obtained from the first one by a reflection.
\end{pf*}

\begin{pf*}{Proof of Proposition~\ref{propnogaus}}
The proof is somewhat similar to that Proposition~\ref{prop<2}.
Suppose that, to the contrary, 
%
\begin{equation}
\label{eqsi} \sigma:=\sqrt{H\bigl(\{0\}\bigr)}>0.
\end{equation}
%
On the other hand, recalling definition (\ref{eqHAB;M}) of $\mathscr{H}
_{p;A,B;M}$ and the conditions $H\in\mathscr{H}_{p;A,B;M}$ and $A>0$,
one sees
that necessarily $\supp H\setminus\{0\}\ne\varnothing$.
So, in view of possible rescaling and reflection, w.l.o.g.
\[
1\in\supp H.
\]

Take now any $\beta\in(0,(\frac{p-2}{p-1})^{1/(p-2)} )$, so that
%
\begin{equation}
\label{eqvp=} \varepsilon:=\frac{\beta^{p-2}}{p-2}\in\biggl(0,\frac
1{p-1}\biggr)\subset(0,1).
\end{equation}
Introduce then
%
\begin{equation}
\label{eqal,ta,tb} {\tilde{a}}:=\frac{1-(p-1)\varepsilon}{(1-\varepsilon
)^p} \quad\mbox{and}\quad {\tilde{b}}:=1-\varepsilon.
\end{equation}
Define the real-valued measure $\Delta:=\Delta_{\beta,\delta}$ by
the condition
%
\begin{eqnarray}\label{eqDeGaus}
\int_\mathbb{R}f(u)\Delta(\dd u)&=&
\frac{1}2 f(\beta)+\frac{1}2 f(-\beta)-f(0)+{\tilde{a}} {\tilde
{b}}f({\tilde{b}})
\nonumber\\[-8pt]\\[-8pt]
&&{} -\frac{1}{H([1-\delta,1+\delta])} \int
_{[1-\delta,1+\delta
]}f(u)H(\dd u)\nonumber
\end{eqnarray}
for all locally bounded (say) Borel functions $f\dvtx \mathbb{R}\to
\mathbb{R}$, where
$\delta$ is any positive real number, so that $H([1-\delta,1+\delta])>0$.
For $\sigma$ as in (\ref{eqsi}), let
\[
t_0:=\sigma^2\wedge H\bigl([1-\delta,1+\delta] \bigr).
\]
Then $t_0>0$ and for all $t\in[0,t_0]$ the measure
%
\begin{equation}
\label{eqHt,be,de} H_t:=H_{\beta,\delta,t}:=H+t\Delta
\end{equation}
is nonnegative.
By Lemma~\ref{lemdiff},
%
\begin{eqnarray}\label{eqLbe}
\qquad \mathscr{L}(\beta) &:=&\lim_{\delta\downarrow0}\lim
_{t\downarrow
0}\frac{\E|X+Y_{H_{\beta,\delta,t}}|^q-\E|X+Y_H|^q}t
\nonumber\\[-8pt]\\[-8pt]
&=&\int_0^1\dd s (1-s) \biggl[
\frac{1}2 h(s\beta)+\frac{1}2 h(-s\beta)-h(0) +{\tilde{a}} {
\tilde{b}}h(s{\tilde{b}})-h(s) \biggr],\nonumber
\end{eqnarray}
where $h$ is still defined by (\ref{eqhq}).
Let now $\beta\downarrow0$.
Then, in view of (\ref{eqal,ta,tb}) and (\ref{eqvp=}), the expression
${\tilde{a}}{\tilde{b}}h(s{\tilde{b}})-h(s)$ in (\ref{eqLbe}) is
$O(\varepsilon)=O(\beta^{p-2})=o(\beta^2)$ uniformly over $s\in[0,1]$.
Concerning the other part of the expression in the brackets in (\ref
{eqLbe}), by~(\ref{eqR==}) and (\ref{eqR=}),
\begin{eqnarray*}
\frac{1}2 h(u)+\frac{1}2 h(-u)-h(0)&=&\frac{u^2}2
\bigl[(R_1h) (0;u)+(R_1h) (0;-u)\bigr]
\\
&=& \frac{u^2}2 \int_{-1}^1
\dd v \bigl(1-|v|\bigr)h''(vu)
\end{eqnarray*}
for all $u\in\mathbb{R}$. So,
\[
\mathscr{L}(\beta)=\frac{\beta^2}2 \int_0^1
\dd s (1-s)s^2\int_{-1}^1\dd v
\bigl(1-|v|\bigr)h''(vs\beta) + o\bigl(\beta^2
\bigr),
\]
whence
\[
\lim_{\beta\downarrow0}\lim_{\delta\downarrow0}\lim
_{t\downarrow0} \frac{\E|X+Y_{H_{\beta,\delta,t}}|^q-\E|X+Y_H|^q}{\beta
^2t/24} 
=h''(0)=
4!
\binom q4 \E|W|^{q-4},
\]
where
$W=X+Y_H$, which is a nondegenerate r.v., so that $\E|W|^{q-4}>0$.

Now (\ref{eqatta}) implies that
eventually
\[
\mathscr{S}_{p,q;\le A,\le B,X;M}=\E|X+Y_H|^q<
\E|X+Y_{H_{\beta,\delta,t}}|^q.
\]
In this context, we say that an assertion $\mathscr{A}=\mathscr
{A}_{\beta,\delta,t}$ holds ``eventually'' if
$\exists\beta_0\in(0,(\frac{p-2}{p-1})^{1/(p-2)} )$
$\forall\beta\in(0,\beta_0)$
$\exists\delta^*_\beta\in(0,\infty)$
$\forall\delta\in(0,\delta^*_\beta)$
$\exists t_{\beta,\delta}\in(0,t_0)$
$\forall t\in(0,t_{\beta,\delta})$ $\mathscr{A}_{\beta,\delta,t}$ holds.

Thus, we obtain a contradiction with the definition of $\mathscr
{S}_{p,q;\le A,\le B,X;M}$ in (\ref{eqS<A<B;M}), because, as we shall
check in moment, $H_{\beta,\delta,t}\in\mathscr{H}_{p;\le A,\le
B;M}$ eventually.
Indeed, by (\ref{eqDeGaus}) and (\ref{eqal,ta,tb}),
\[
\int_\mathbb{R}\Delta(\dd x)={\tilde{a}} {
\tilde{b}}-1=\frac
{1-(p-1)\varepsilon}{(1-\varepsilon)^{p-1}}-1<0,
\]
so that 
\[
\int_\mathbb{R}H_{\beta,\delta,t}(\dd x)=\int
_\mathbb{R}H(\dd x)+t\int_\mathbb{R}\Delta(\dd
x)<\int_\mathbb{R}H(\dd x)\le B,
\]
by (\ref{eqHt,be,de}) and (\ref{eqH<A<B;M}).
Next, by (\ref{eqal,ta,tb}),
\begin{eqnarray*}
&& \lim_{\delta\downarrow0}\int_\mathbb{R}|x|^{p-2}
\Delta(\dd x)
\\
&&\qquad = \beta^{p-2}+{\tilde{a}} {\tilde{b}}^{p-1}-1
= (p-2)\varepsilon+\frac{1-(p-1)\varepsilon}{1-\varepsilon}-1
\\
&&\qquad =-\frac{(p-2)\varepsilon^2}{1-\varepsilon}<0.
\end{eqnarray*}
It follows that eventually
\begin{eqnarray*}
&& \int_\mathbb{R}|x|^{p-2}H_{\beta,\delta,t}(\dd x)
\\
&&\qquad =\int
_\mathbb{R}|x|^{p-2}H(\dd x)+t\int
_\mathbb{R}|x|^{p-2}\Delta(\dd x) <\int
_\mathbb{R}|x|^{p-2}H(\dd x)\le A,
\end{eqnarray*}
again by (\ref{eqHt,be,de}) and (\ref{eqH<A<B;M}).

Also, the conditions $H\in\mathscr{H}_{p;\le A,\le B;M}$ and $1\in
\supp H$
imply $M\ge1$.
So, eventually $\supp H_{\beta,\delta,t}\subseteq\supp H\cup\{\beta,-\beta,{\tilde{b}}\}\subseteq[-M,M]$ in view of (\ref
{eqal,ta,tb}). 

By (\ref{eqH<A<B;M}), we conclude that indeed $H_{\beta,\delta,t}\in
\mathscr{H}_{p;\le A,\le B;M}$ eventually.
Thus, assumption (\ref{eqsi}) 
leads to a contradiction.
\end{pf*}

\begin{pf*}{Proof of Proposition~\ref{propcard=1}}
By Propositions~\ref{prop<2} and \ref{propnogaus} and the condition %
$\supp H\subseteq(-M,M)$,
%
\begin{equation}
\label{eqH*=} H=w_1\delta_{c_1}+w_2
\delta_{-c_2}
\end{equation}
for some $c_1$ and $c_2$ in the interval $(0,M)$ and some nonnegative
real $w_1$ and $w_2$ such that
%
\begin{equation}
\label{eq=B,=A} w_1+w_2=B\quad\mbox{and}\quad
c_1^{p-2}w_1+c_2^{p-2}w_2=A.
\end{equation}

It is enough to show that $w_1\wedge w_2=0$.
To obtain a contradiction, suppose the contrary,
%
\begin{equation}
\label{eqw} w:=w_1\wedge w_2>0.
\end{equation}
Then, by the implicit function theorem, there exist a real number $\tau
_*>0$ and an infinitely differentiable mapping $(-\tau_*,\tau_*)\ni
\tau\mapsto(\tc_1(\tau),\tc_2(\tau))$
such that
%
\begin{equation}
\label{eqtc0} \tc_1(0)=c_1,\qquad\tc_2(0)=c_2,
\end{equation}
and
for each $\tau\in(-\tau_*,\tau_*)$ one has $\tc'_1(\tau)\tc
'_2(\tau)\ne0$,
%
\begin{equation}
\label{eqtilde} 0<\tc_1(\tau),\tc_2(\tau)<M\quad
\mbox{and}\quad\tc_1(\tau)^{p-2}+\tc_2(
\tau)^{p-2}=c_1^{p-2}+c_2^{p-2}.
\end{equation}
(In this case, this mapping could also be defined explicitly, e.g., by
the formulas $\tc_1(\tau)=(c_1^{p-2}+\tau)^{1/(p-2)}$ and $\tc
_2(\tau)=(c_2^{p-2}-\tau)^{1/(p-2)}$, with $\tau_*=\frac{1}2\min
[M^{p-2}-c_1^{p-2},M^{p-2}-c_2^{p-2},c_1^{p-2},c_2^{p-2}]$.)
Note that the\vspace*{1pt} condition\break $\tc'_1(\tau)\tc'_2(\tau)\ne0$, taken
together with (\ref{eqtilde}), implies
%
\begin{equation}
\label{eqprod<0} \tc'_1(\tau)\tc'_2(
\tau)<0.
\end{equation}
By choosing a possibly smaller real $\tau_*>0$, let us assume w.l.o.g.
that, on the interval $(-\tau_*,\tau_*)$, the derivatives of any
order of the functions $\tc_1$ and $\tc_2$ are each uniformly
continuous and hence bounded, and also that the functions $\tc_1$ and
$\tc_2$ are each positive and bounded away from $0$.

For each $\tau\in(-\tau_*,\tau_*)$, introduce the real-valued measure
%
\begin{equation}
\label{eqDetau} \Delta_\tau:=\delta_{\tc_1(\tau)}+
\delta_{-\tc_2(\tau)}-\delta_{c_1}-\delta_{-c_2}
\end{equation}
and then the measures
%
\begin{equation}
\label{eqHt,tau} H_{t,\tau}
:=H+t\Delta_\tau\qquad
\mbox{for } t\in(-w,w),
\end{equation}
where $w$ is as in (\ref{eqw}).
By (\ref{eqHAB;M}), (\ref{eqHAB}), (\ref{eqH*=}), (\ref{eq=B,=A})
and (\ref{eqtilde}), these measures are all in~$\mathscr{H}_{p;A,B;M}$.

In the rest of this proof, it is assumed that $\tau\in(-\tau_*,\tau
_*)$, $t\in(-w,w)$, $\{j,k,\ell\}\subset\{1,2,3,4\}$, and $x\in
\mathbb{R}$---unless otherwise indicated.

Letting now
%
\begin{equation}
\label{eqgt,tau} g_{t,\tau}(x):=\E|x+X+Y_{H_{t,\tau}}|^q,
\end{equation}
then using Lemma~\ref{lemdiff} and recalling (\ref{eqDetau}), one has
%
\begin{eqnarray}\label{eqD2g}
\mathscr{D}(\tau) &:=& \pd{{}^2g_{t,\tau}(0)}
{t^2} \bigg|_{t=0}\nonumber
\\
&=& \int
_\Omega\mu(\dd\omega) \int_{\mathbb{R}^2}
\Delta_\tau(\dd u_1)\Delta_\tau(\dd
u_2) h(x+su_1+tu_2)
\\
&=& \int_\Omega\mu(\dd\omega) F_\omega(\tau),\nonumber
\end{eqnarray}
where
%
\begin{eqnarray}
\Omega&:=&(0,1)^2\times\mathbb{R},\qquad\omega:=(s,t,x)\in\Omega,
\nonumber
\\
\mu(\dd\omega)&:=&\dd{t} \,\dd{s} (1-t) (1-s)\P(X+Y_H\in\dd x),\label{eqmu}
\\
F_\omega(\tau)&:=& \sum_{j,k=1}^4 v_j v_k h \bigl(x+sb_j(\tau)+tb_k(\tau) \bigr), \label{eqF=}
\\
h(x)&:=&  4!\binom q4 |x|^{q-4}, \label{eqhp}
\\
(v_1,v_2,v_3,v_4)&:=&(1,1,-1,-1),\label{eqvv}
\\
\bigl(b_1(\tau),b_2(\tau),b_3(
\tau),b_4(\tau) \bigr)&:=& \bigl(\tc_1(\tau),-
\tc_2(\tau),c_1,-c_2 \bigr). \label{eqbb}
\end{eqnarray}
Next,
%
\begin{equation}
\label{eqM39108} \mathscr{D}(\tau)=\sum_{j,k=1}^4
v_j v_k \mathscr{D}_{j,k}(\tau),
\end{equation}
where
%
\begin{equation}
\label{eqMjkx} \mathscr{D}_{j,k}(\tau):=\int_\Omega
\mu(\dd\omega)h \bigl(x+sb_j(\tau)+tb_k(\tau) \bigr).
\end{equation}
By Lemma~\ref{lemunderint,simpl},
%
\begin{equation}
\label{eqMjk=} \mathscr{D}'_{j,k}(\tau) =\int
_\Omega\mu(\dd\omega) f(\omega,\tau), 
\end{equation}
where
%
\begin{equation}
\label{eqfom,tau} f(\omega,\tau):=f_{j,k}(\omega,\tau):=h'
\bigl(x+sb_j(\tau)+tb_k(\tau) \bigr)
\bigl[sb'_j(\tau)+tb'_k(\tau)
\bigr],
\end{equation}
which is clearly bounded in $(\omega,\tau)\in\Omega\times(-\tau
_*,\tau_*)$.
For each $\varepsilon\in[0,\infty)$, introduce the set
%
\begin{equation}
\label{eqOmvp} \Omega_\varepsilon:=\Omega_{j,k;\varepsilon}:= \bigl\{
\omega=(s,t,x)\in\Omega\dvtx \bigl|x+sb_j(0)+tb_k(0)\bigr|>
\varepsilon\bigr\}.
\end{equation}
Since $b_j(\tau)$ is uniformly continuous in $\tau\in(-\tau_*,\tau
_*)$ for each $j$, one sees that $|x+sb_j(\tau)+tb_k(\tau)|$ is
continuous in $\tau$ uniformly over all $(\omega,\tau,j,k)\in\Omega
\times(-\tau_*,\tau_*)\times\{1,2,3,4\}\times\{1,2,3,4\}$.
So, by further decreasing (if necessary) the value of $\tau_*>0$, let
us assume, again w.l.o.g., that
%
\begin{eqnarray}
\qquad \bigl|x+sb_j(\tau)+tb_k(\tau)\bigr|&>&\varepsilon/2\qquad\mbox{for
all } (\omega,\tau)\in\Omega_\varepsilon\times(-\tau_*,\tau_*),
\label{eq>vp2}
\\
\bigl|x+sb_j(\tau)+tb_k(\tau)\bigr|&\le&2\varepsilon\qquad
\mbox{for all } (\omega,\tau)\in(\Omega\setminus\Omega_\varepsilon
)\times(-
\tau_*,\tau_*). \label{eq<2vp}
\end{eqnarray}
By (\ref{eqfom,tau}), for $(\omega,\tau)\in\Omega_0\times(-\tau
_*,\tau_*)$, the partial derivative of $f(\omega,\tau)$ in $\tau$ is
%
\begin{equation}
\label{eqp2f=} (\partial_2f) (\omega,\tau)=D_1(\omega,
\tau)+D_2(\omega,\tau),
\end{equation}
where
\begin{eqnarray*}
D_1(\omega,\tau)&:=&h' \bigl(x+sb_j(
\tau)+tb_k(\tau) \bigr) \bigl[sb''_j(
\tau)+tb''_k(\tau)\bigr],
\\
D_2(\omega,\tau)&:=&h''
\bigl(x+sb_j(\tau)+tb_k(\tau) \bigr)
\bigl[sb'_j(\tau)+tb'_k(\tau)
\bigr]^2.
\end{eqnarray*}
In view of the condition $q>5$, definition (\ref{eqhp}), inequality
(\ref{eq>vp2}), and the boundedness of all the derivatives of the
functions $b_j$ on the interval $(-\tau_*,\tau_*)$,
%
\begin{eqnarray}\label{eqD1,D2<}
\lleft. %
\begin{array} {l} \bigl|D_1(\omega,
\tau)\bigr| \le K \bigl(1+|x|^{q-5} \bigr) \quad\mbox{and}
\\[4pt]
\bigl|D_2(\omega,\tau)\bigr|
\le K \bigl(1+|x|^{(q-6)_+}+\varepsilon
^{-(6-q)_+} \bigr)
\end{array}
\rright\}
\nonumber\\[-8pt]\\[-8pt]
\eqntext{\mbox{for all }(\omega,\tau)\in\Omega_\varepsilon\times(-\tau_*,\tau_*);}
\end{eqnarray}
here and in the rest of this proof, $K$ denotes various positive real
constants which do not depend on $\omega$, $\tau$, or $\varepsilon$.
So, by (\ref{eqp2f=}),
%
\begin{eqnarray}\label{eqintD1to0}
\bigl|(\partial_2f) (\omega,\tau)\bigr|\le g_\varepsilon(
\omega):=K\bigl(1+|x|^{q-5}+\varepsilon^{-(6-q)_+}\bigr)
\nonumber\\[-8pt]\\[-8pt]
\eqntext{\mbox{for all } (\omega,\tau)\in\Omega_\varepsilon\times(-\tau_*,\tau_*).}
\end{eqnarray}
%
By (\ref{eqHt,tau}), (\ref{eqH*=}), (\ref{eqDetau}), (\ref{eqcosh}),
and (\ref{eqX}),
$\int_{\Omega_\varepsilon}|\dd\mu| g_\varepsilon<\infty$, where
$\mu$ is still as in (\ref{eqmu}).

Next, 
by (\ref{eqD1,D2<}) and dominated convergence,
%
\begin{equation}
\label{eqintD1} \sup_{\tau\in(-\tau_*,\tau_*)}\int_{\Omega\setminus
\Omega_\varepsilon}\bigl|\mu(\dd\omega)D_1(\omega,\tau)\bigr|
\underset{\varepsilon
\downarrow0} {\longrightarrow}0.
\end{equation}
Further, $|D_2(\omega,\tau)|\le K |x+sb_j(\tau)+tb_k(\tau)|^{q-6}$
for all $(\omega,\tau)\in\Omega_0\times(-\tau_*,\tau_*)$,
whence, by (\ref{eq<2vp}),
with $\nu(\dd x):=\P(X+Y_H\in\dd x)$,
\begin{eqnarray*}
&& \int_{\Omega\setminus\Omega_\varepsilon}\bigl| \mu(\dd\omega)D_2(\omega,\tau)\bigr|
\\
&&\qquad \le K\int_0^1\dd t
\int_\mathbb{R}\nu(\dd x)
\int_0^1\dd s \bigl|x+sb_j(
\tau)+tb_k(\tau)\bigr|^{q-6}
\\
&&\hspace*{124pt}{}\times  \mathrm{I}\bigl\{\bigl|x+sb_j(
\tau)+tb_k(\tau)\bigr|\le2\varepsilon\bigr\}
\\
&&\qquad \le K\int_0^1\dd t \int_\mathbb{R}
\nu(\dd x) 
\frac{1}{b_j(\tau)} \int_\mathbb{R}\dd v
|v|^{q-6} \mathrm{I}\bigl\{ |v|\le2\varepsilon\bigr\} =\frac{2K(2\varepsilon
)^{q-5}}{b_j(\tau)(q-5)}
\underset{\varepsilon\downarrow0} {\longrightarrow}0
\end{eqnarray*}
uniformly in $\tau\in(-\tau_*,\tau_*)$, since the functions $b_j$
are bounded away from $0$ on $(-\tau_*,\tau_*)$.
Combining this with (\ref{eqp2f=}) and (\ref{eqintD1}), one has
\[
\sup_{\tau\in(-\tau_*,\tau_*)}\int_{\Omega\setminus\Omega
_\varepsilon}\bigl|\mu(\dd\omega) (
\partial_2f) (\omega,\tau)\bigr| \underset{\varepsilon\downarrow0} {
\longrightarrow}0.
\]
Therefore and by (\ref{eqintD1to0}), one may use Lemma~\ref{lemunderint}
together with (\ref{eqMjk=}) and (\ref{eqfom,tau})
to conclude that
$\mathscr{D}''_{j,k}(0)
=
\int_\Omega\mu(\dd\omega) \pd{{}^2}{\tau^2} h (x+sb_j(\tau
)+tb_k(\tau) ) |_{\tau=0}$ and\vspace*{2pt} hence, by~(\ref{eqM39108}),
(\ref{eqMjkx}) and (\ref{eqF=}), 
%
\begin{equation}
\label{eqM399} \mathscr{D}''(0) = \int
_\Omega\mu(\dd\omega) F''_\omega(0);
\end{equation}
that is, we have shown that the second integral expression of $\mathscr
{D}(\tau)$ in (\ref{eqD2g}) can be twice differentiated (at least at
$\tau=0$) under the integral sign to obtain the corresponding integral
expression of $\mathscr{D}''(0)$.
Note here that $F''_\omega(0)$ is defined only for 
$\omega\in\bigcap_{j,k=1}^4\Omega_{j,k;0}
$, where $\Omega_{j,k;0}$ is understood according to
(\ref{eqOmvp}).
However, this causes no problem, since $\mu(\Omega\setminus
\bigcap_{j,k=1}^4\Omega_{j,k;0}
)=0$.\vspace*{1pt} 

In view of (\ref{eqF=}), (\ref{eqhp}), (\ref{eqvv}), (\ref{eqbb}),
and (\ref{eqtc0}), it is straightforward but tedious to check that %
%
\begin{eqnarray}
\qquad F_\omega(0)&=&0, \label{eqF0=0}
\\
F'_\omega(0)&=&0, \label{eqF1}
\\
F''_\omega(0)&=&2st \bigl
\{h'' \bigl(x-(s+t)c_2 \bigr)
\tc'_2(0)^2 +h''
\bigl(x+(s+t)c_1 \bigr)\tc'_1(0)^2
\nonumber\\[-8pt] \label{eqF2}\\[-8pt]
&&\hspace*{17pt}{} - \bigl[h''(x+sc_1-tc_2)+h''(x-sc_2+tc_1)
\bigr] \tc'_1(0)\tc'_2(0)
\bigr\}\nonumber
\end{eqnarray}
for all 
$\omega\in\bigcap_{j,k=1}^4\Omega_{j,k;0}$.
The equality in (\ref{eqF1}) in fact holds for all $\omega\in\Omega$
and any continuously differentiable function $h$, not necessarily the
one defined by (\ref{eqhp}),
whereas the equality in (\ref{eqF0=0}) holds for any function $h\dvtx
\mathbb{R}\to\mathbb{R}$ whatsoever.

By (\ref{eqhp}), $h''(z)>0$ for all real $z\ne0$. So, by (\ref{eqF2})
and (\ref{eqprod<0}), $F''_\omega(0)>0$ for all 
$\omega\in\bigcap_{j,k=1}^4\Omega_{j,k;0}$.
It follows 
by
(\ref{eqD2g}), (\ref{eqF0=0}), (\ref{eqF1}), and Lemma~\ref{lemunderint,simpl} that $\mathscr{D}(0)=\mathscr{D}'(0)=0$,
whereas, by (\ref{eqM399}), 
$\mathscr{D}''(0)>0$
and hence
$\mathscr{D}(\tau)>0$ for some $\tau\in(-\tau_*,\tau_*)$ (in fact
for all nonzero $\tau$ close enough to $0$).
Take any such $\tau$. Then, by (\ref{eqD2g}),
\[
\pd{{}^2g_{t,\tau}(0)} {t^2} \bigg|_{t=0}=
\mathscr{D}(\tau)>0,
\]
which implies that $g_{0,\tau}(0)<g_{-t,\tau}(0)\vee g_{t,\tau}(0)$
if $|t|$ is small enough. In view of~(\ref{eqHt,tau}) and (\ref
{eqgt,tau}), this means that for 
all $t\in(-w,w)$ with small enough $|t|$,
\[
\E|X+Y_H|^p=\E|X+Y_{H_{0,\tau}}|^p<
\E|X+Y_{H_{-t,\tau}}|^p\vee\E|X+Y_{H_{t,\tau}}|^p,
\]
which is a contradiction, in view of the conditions $H\in\mathscr{H}
_{*,p,q;A,B;X;M}$ and $H_{t,\tau}\in\mathscr{H}_{p;A,B;M}$ for all
$(\tau,t)\in(-\tau_*,\tau_*)\times(-w,w)$, and the definition (\ref
{eqH*}) of $\mathscr{H}_{*,p,q;A,B;X;M}$.
\end{pf*}

\begin{pf*}{Proof of Proposition~\ref{proplim}}
Since $|c_2|\to\infty$, w.l.o.g. $c_2\ne0$.
So, the definition $\lambda_2:=w_2/c_2^2$ makes sense, and $\lambda
_2\in[0,\infty)$. If $\lambda_2=0$, let $\tPi_{\lambda_2}:=0$.
So, $Y_H\stackrel{\mathrm{D}}=Y_{w_1\delta_{c_1}}+c_2\tPi_{\lambda
_2}$ and hence
%
\begin{equation}
\label{eq=sumTj} e^{\lambda_2}\E|X+Y_H|^q=\sum
_{j=0}^\infty T_j,
\end{equation}
where
\[
T_j:=\frac{\lambda_2^j}{j!} \E\bigl|X+Y_{w_1\delta_{c_1}}+c_2(j-
\lambda_2)\bigr|^q,
\]
letting $\lambda_2^0:=1$ even if $\lambda_2=0$.
So, in view of the conditions $|c_2|\to\infty$, $|c_2|^{q-2}w_2\to
a$, $w_1+w_2=B$, and $c_1\to b$, one has
$w_2\to0$, $w_1\to B$,
$\lambda_2\to0$, $Y_{w_1\delta_{c_1}}\stackrel{\mathrm
{D}}\longrightarrow Y_{B\delta_b}$,
$c_2\lambda_2\to0$, $|c_2|^q\lambda_2\to a$, $|c_2|^q\lambda_2^2\to
0$, whence,
by dominated convergence,
%
\begin{eqnarray}
\label{eqT0,T1} T_0&\to&\E|X+Y_{B\delta_b}|^q,
\nonumber\\[-8pt]\\[-8pt]
T_1&=&\E\bigl|\lambda_2^{1/q}(X+Y_{w_1\delta_{c_1}})+
\bigl(|c_2|^q\lambda_2\bigr)^{1/q}
\sign c_2 \,(1-\lambda_2) \bigr|^q\to a.\nonumber
\end{eqnarray}
Also, eventually $\lambda_2\in[0,1]$ and hence
%
\begin{equation}
\label{eqotherTj} 2^{1-q}\sum_{j=2}^\infty|T_j|
\le\lambda_2^2\sum_{j=2}^\infty
\frac{1}{j!} \E|X+Y_{w_1\delta_{c_1}}|^q +|c_2|^q
\lambda_2^2\sum_{j=2}^\infty
\frac{j^q}{j!}\longrightarrow0.
\end{equation}
Combining (\ref{eq=sumTj}), (\ref{eqT0,T1}), (\ref{eqotherTj}) and
recalling that $\lambda_2\to0$, one
completes the proof.
\end{pf*}

\subsection{Conclusion of the proof of Theorem~\texorpdfstring{\protect\ref{th}}{1.3}}\label{conclproof}
Consider first the case when $X$ is bounded and
%
\begin{equation}
\label{eq5<q<p} 5<q<p.
\end{equation}

Recall definition (\ref{eqH*}) of $\mathscr{H}_{*,p,q;A,B;X;M}$.
By Propositions~\ref{propcard=1} and \ref{propnogaus}, for each real
$M$ as in (\ref{eqM>}), either
$
\mathscr{H}_{*,p,q;A,B;X;M}\subseteq\{B\delta_c,B\delta_{-c}\}$ or
there is
some $H_{*,M}\in\mathscr{H}_{*,p,q;A,B;X;M}$ such that $\supp
H_{*,M}=\{
-c_{2,M},c_{1,M}\}$ for some real $c_{1,M}$ and $c_{2,M}$ such that
$0<c_{1,M}\wedge c_{2,M}\le c_{1,M}\vee c_{2,M}=M$.
So, w.l.o.g. one of the following two cases holds:
\begin{longlist}[\textit{Case}~1.]
\item[\textit{Case}~1.] $
\mathscr{H}_{*,p,q;A,B;X;M}\subseteq\{B\delta_c,B\delta_{-c}\}$ for
all real
$M\ge c$.
\item[\textit{Case}~2.] There exist 
sequences $(M_k)$ in $[c,\infty)$, $(b_k)$ in $[0,c]$, $(w_{1,k})$ in
$[0,B]$, and $(w_{2,k})$ in $[0,B]$ such that $M_k\uparrow\infty$, and
for all $k$ one has $H_k:=w_{1,k}\delta_{M_k}+w_{2,k}\delta_{-b_k}\in
\mathscr{H}_{*,p,q;A,B;X;M_k}$, $w_{1,k}+w_{2,k}=B$, and
$b_k^{p-2}w_{1,k}+M_k^{p-2}w_{2,k}=A$.
\end{longlist}

In case~1, by (\ref{eqH*}),
%
\begin{eqnarray}\label{eqS=max}
\mathscr{S}_{p,q;A,B;X;M}=\max\bigl(\E|X+Y_{B\delta_c}|^q,
\E|X+Y_{B\delta_{-c}}|^q \bigr)
\nonumber\\[-8pt]\\[-8pt]
\eqntext{\mbox{for all real }M\ge c.}
\end{eqnarray}

Let us show that (\ref{eqS=max}) holds in case~2 as well.
W.l.o.g.,
$b_k\to b$ for some 
$b\in[0,c]$.
Also, $0\le M_k^{q-2}w_{2,k}\le M_k^{q-2}\frac{A}{M_k^{p-2}}\to0$, by
the condition $q<p$ in (\ref{eq5<q<p}). So, by Proposition~\ref{proplim},
%
\begin{equation}
\label{eqcase2} \E|X+Y_{H_k}|^q\to\E|X+Y_{B\delta_{-b}}|^q.
\end{equation}
Since $H_k\in\mathscr{H}_{*,p,q;A,B;X;M_k}$ and $\mathscr{S}_{p,q;A,B;X;M}$
is obviously nondecreasing in $M>0$, it now follows that
%
\begin{eqnarray}\label{eqsuple}
\sup_{M>0}\mathscr{S}_{p,q;A,B;X;M} &=&
\E|X+Y_{B\delta_{-b}}|^q
\nonumber\\[-8pt]\\[-8pt]
&\le&\mathscr{S}_{p,q;Bb^{p-2},B;X;c} \le
\mathscr{S}_{p,q;A,B;X;c}.\nonumber
\end{eqnarray}
The last inequality follows by Proposition~\ref{propincr}, because
$b\in[0,c]$ and hence $Bb^{p-2}\le Bc^{p-2}=A$.
Moreover, if $b\in[0,c)$ then, again by Proposition~\ref{propincr},
the last inequality in (\ref{eqsuple}) is strict, which is a contradiction.
Thus, necessarily $b=c$, and so, by the equality in (\ref{eqsuple}),
(\ref{eqS=max}) holds in case~2 as well, because obviously $\mathscr
{S}_{p,q;A,B;X;M}\ge\E|X+Y_{B\delta_c}|^q$ for all real $M\ge c$.

Take now any real $M\ge c$ and any $H\in\mathscr{H}_{p;\le A,\le
B;M}$. Then,
by (\ref{eqS=max}),
(\ref{eqle=eq}), and~(\ref{eqS<A<B;M}),
%
\begin{equation}
\label{eqqin[5,p]} \E|X+Y_H|^q\le\max\bigl(
\E|X+Y_{B\delta_c}|^q,\E|X+Y_{B\delta
_{-c}}|^q \bigr)
\end{equation}
---provided that $q\in(5,p)$.
Since $\E|X+Y_H|^q$ is continuous in $q\in(0,\infty)$ [cf. the
second equality in (\ref{eqattain})], inequality (\ref{eqqin[5,p]})
holds for all $q\in[5,p]$---provided that $p>5$.

Let us show that (\ref{eqqin[5,p]}) holds when $p=5$ (and then $q=5$
as well).
%
Take any $H\in\mathscr{H}_{5,\le A,\le B;M}$ and any sequence $(p_n)$ in
$(5,\infty)$ such that $p_n\downarrow5$ as $n\to\infty$.
Then
$|x|^{p_n-2}\to|x|^{5-2}$ uniformly in $x\in[-M,M]$ and hence
\[
A_n:=A\vee\int_\mathbb{R}|x|^{p_n-2}H(\dd
x)\quad\longrightarrow\quad A\vee\int_\mathbb{R}|x|^{5-2}H(\dd
x)=A.
\]
So, recalling (\ref{eqla,c}) and letting $b_n:=c_{p_n}(A_n,B)$, 
one has $b_n\to c$.
Also, clearly $H\in\mathscr{H}_{p_n;\le A_n,\le B;M}$ for all $n$. Therefore,
by 
(\ref{eqqin[5,p]}) with $q=5$, 
%
\begin{eqnarray}
\label{eq5} \E|X+Y_H|^5
&\le& \max\bigl(
\E|X+Y_{B\delta_{b_n}}|^5,\E|X+Y_{B\delta
_{-b_n}}|^5
\bigr) 
\nonumber\\[-8pt]\\[-8pt]
&\to&\max\bigl(\E|X+Y_{B\delta_c}|^5,
\E|X+Y_{B\delta_{-c}}|^5 \bigr).\nonumber
\end{eqnarray}
Thus, indeed (\ref{eqqin[5,p]}) holds when $p=q=5$.

Take now any ${\mathbf{X}}=(X_1,\ldots,X_n)\in\mathscr{X}_{p;X;\le
A,\le B}$
and abandon the assumption that the r.v. $X$ is bounded.
Let $X_0:=X$.
By Proposition~\ref{proptrunc}, for each $i\in\{0,\ldots,n\}$ and
each real $M>0$ there is a truncated version $X_{i,M}$ of $X_i$ such that:
\begin{longlist}[(iii)]
\item[(i)] $\E X_{i,M}=0$;
\item[(ii)] $|X_{i,M}|\le M\wedge|X_i|$;
\item[(iii)] $\E f(X_{i,M})\le\E f(X_i)$ for all convex
functions $f\dvtx \mathbb{R}\to\mathbb{R}$;\label{jensen}
\item[(iv)]$X_{i,M}\to X_i$ a.s. as $M\to\infty$;
\item[(v)]$X_{0,M},\ldots,X_{n,M}$ are independent.
\end{longlist}
Then obviously
%
\begin{equation}
\label{eqleA,B} (X_{1,M},\ldots,X_{n,M})\in
\mathscr{X}_{p;X;\le A,\le B}.
\end{equation}
%
Letting now $S_M:=X_{1,M}+\cdots+X_{n,M}$,
one also has
$|X_{0,M}+S_M|^q\le
(n+1)^{q-1} (|X_{0,M}|^q+\sum_1^n|X_{i,M}|^q )\le
(n+1)^{q-1} (|X|^q+\sum_1^n|X_i|^q )$.
%
So, by dominated convergence,
%
\begin{equation}
\label{eqconv} \E|X_{0,M}+S_M|^q\underset{M
\to\infty} {\longrightarrow}\E|X+S_{\mathbf{X}}|^q.
\end{equation}


On the other hand, by Theorem~\ref{teoA} (with $\E|X_{0,M}+\cdot\,|^p$ and $X_{i,M}$ in place of $f$ and $X_i$)
and (\ref{eqcf}),
%
\begin{equation}
\label{eqS<Y} \E|X_{0,M}+S_M|^q\le
\E|X_{0,M}+Y_{H_{*,M}}|^q,
\end{equation}
where
\[
H_{*,M}(E):=\int_E x^2\sum
_1^n\P(X_{i,M}\in\dd x)
\]
for all Borel sets $E\subseteq\mathbb{R}$.
It follows from (\ref{eqleA,B}) that the measure $H_{*,M}$ is in
$\mathscr{H}
_{p;\le A,\le B;M}$.
By (\ref{eqS<Y}), (\ref{eqqin[5,p]}) (proved for bounded $X$ and
$H\in\mathscr{H}_{p;\le A,\le B;M}$) and item (iii) on
page~\pageref{jensen},
\begin{eqnarray*}
\E|X_{0,M}+S_M|^q
&\le&
\E|X_{0,M}+Y_{H_{*,M}}|^q
\\
&\le& \max\bigl(
\E|X_{0,M}+Y_{B\delta_c}|^q,\E|X_{0,M}+Y_{B\delta
_{-c}}|^q
\bigr)
\\
&\le&\max\bigl(\E|X+Y_{B\delta_c}|^q,\E|X+Y_{B\delta_{-c}}|^q
\bigr)
\\
&=&\max\bigl(\E|X+c\tPi_{\lambda}|^q,\E|X-c
\tPi_{\lambda}|^q \bigr),
\end{eqnarray*}
where again $\lambda$ and $c$ are as in (\ref{eqla,c}).

Now (\ref{eqconv}) yields
%
\begin{equation}
\label{eq<rhs} \E|X+S_{\mathbf{X}}|^q\le\max\bigl(\E|X+c
\tPi_{\lambda}|^q,\E|X-c\tPi_{\lambda}|^q
\bigr).
\end{equation}
Thus, the first supremum in (\ref{eq}) is no greater than the
right-hand side of (\ref{eq<rhs}).

To complete the proof of Theorem~\ref{th}, it remains to note
that the second supremum in (\ref{eq}) is no less than
the right-hand side of (\ref{eq<rhs}).
Indeed, by
Lemma~\ref{lemlim} with $G=\lambda\delta_c$, one has a sequence
$(\ZZ_n)$ in $\mathscr{X}_{p;A,B}$ such that $S_{\ZZ_n}\stackrel
{\mathrm{D}}\longrightarrow c\tPi
_{\lambda}$.
Now, by the Fatou lemma for the convergence in distribution (Theorem~5.3 in \cite{billingsley}),
$
\liminf_n\E|X+S_{\ZZ_n}|^q\ge\E|X+c\tPi_{\lambda}|^q$, 
so that the second supremum in (\ref{eq}) is no less than $\E|X+c\tPi
_{\lambda}|^q$.
Quite similarly, that supremum is no less than $\E|X-c\tPi_{\lambda
}|^q$, and thus it is indeed no less than
the right-hand side of (\ref{eq<rhs}).

\section{Other proofs} \label{proofs,props}
\mbox{}
\begin{pf*}{Proof of Proposition~\ref{propEE}}
That $\varnothing\ne\mathscr{X}_{p;A,B}$ is part of Lemma~\ref{lemlim},
and the inclusion $\mathscr{X}_{p;A,B}\subseteq\mathscr{X}_{p;\le
A,\le B}$ is trivial.
The homogeneity property holds because for any ${\mathbf{X}}\in
\mathscr{X}_{p;A,B}$ and any real $\kappa>0$, one has $\kappa{\mathbf{X}}\in
\mathscr{X}_{p;\kappa^p A,\kappa^2 B}$.

Now it follows easily by Jensen's inequality 
that $\mathscr{E}_{p;A,B}$
is nondecreasing in $A$ and in $B$.
Indeed, let us first take any $\tA\in(0,A)$ and $\tB\in(0,B)$. Take
then any independent finite sequences ${\mathbf{X}}=(X_1,\ldots,X_n)\in
\mathscr{X}_{p;\tA,\tB}$ and ${\mathbf{Y}}=(Y_1,\ldots,Y_m)\in\mathscr{X}_{p;A-\tA,B-\tB}$; by the already verified first sentence of
Proposition~\ref{propEE}, such ${\mathbf{X}}$ and~${\mathbf{Y}}$ exist.
Then $\ZZ:=(X_1,\ldots,X_n,Y_1,\ldots,Y_m)\in\mathscr{X}_{p;A,B}$.
Moreover, by Jensen's inequality, $\E|S_{\mathbf{X}}|^p\le\E
|S_{\mathbf{X}}+S_{\mathbf{Y}}
|^p=\E|S_\ZZ|^p$.
Thus, $\mathscr{E}_{p;\tA,\tB}\le\mathscr{E}_{p;A,B}$, for any $\tA
\in(0,A)$ and $\tB\in(0,B)$.

This and the homogeneity property in turn imply that $\mathscr
{E}_{p;A,\tB}\le\mathscr{E}_{p;\kappa^p A,\kappa^2 B}
=\kappa^p\mathscr{E}_{p;A,B}$ for any $\tB\in(0,B]$ and any real
$\kappa
>1$. Letting now $\kappa\downarrow1$ and recalling that, by (\ref
{eqros}), $\mathscr{E}_{p;A,B}<\infty$, one concludes that $\mathscr
{E}_{p;A,\tB}\le\mathscr{E}_{p;A,B}$ for any $\tB\in(0,B]$.
Similarly, $\mathscr{E}_{p;\tA,B}\le\mathscr{E}_{p;A,B}$ for\vspace*{2pt} any
$\tA\in(0,A]$.
Thus, indeed $\mathscr{E}_{p;A,B}$ is nondecreasing in $A$ and in $B$.
Now (\ref{eqEEle=}) immediately follows.
\end{pf*}

\begin{pf*}{Proof of Proposition~\ref{propEE,C}}
For brevity, let $K_{A,B}:=\max(\gamma A,B^{p/2})^{1/p}$. Then %
$A/K_{A,B}^p\le1/\gamma$, $B/K_{A,B}^2\le1$, and, in view of (\ref{eqga})
and the homogeneity and monotonicity properties of $\mathscr
{E}_{p;A,B}$ presented in Proposition~\ref{propEE},
\[
C_{p;\gamma}=\sup_{A,B>0}K_{A,B}^{-p}
\mathscr{E}_{p;A,B} =\sup_{A,B>0}\mathscr{E}_{p;A/K_{A,B}^p,B/K_{A,B}^2}
\le\mathscr{E}_{p;1/\gamma,1}.
\]
%
On the other hand, by (\ref{eqga}), $\mathscr{E}_{p;1/\gamma,1}\le
C_{p;\gamma}$.
Thus, the first equality in Proposition~\ref{propEE,C} is verified.

The second equality there easily follows from (and in fact is
equivalent to) the first one.
Indeed,\vspace*{1pt} choosing $\gamma=B^{p/2}/A$ and using again the homogeneity
property, one has
$
\mathscr{E}_{p;A,B}=B^{p/2}\mathscr{E}_{p;1/\gamma,1}=B^{p/2}C_{p;\gamma}
=B^{p/2}C_{p;B^{p/2}/A}$. 
\end{pf*}

\begin{pf*}{Proof of Theorem~\ref{th2<p<3}}
Take any ${\mathbf{X}}\in\mathscr{X}_{p;X;\le A,\le B}$.
Let $\sigma:=\sqrt{\Var S_{\mathbf{X}}}$, so that $\sigma\in
[0,\sqrt B ]$. If $\sigma
=0$ then, by Jensen's inequality,
$\E|X+S_{\mathbf{X}}|^p=\E|X|^p\le\E|X+B^{1/2}Z|^p
\le A+\E|X+B^{1/2}Z|^p$, whence
%
\begin{equation}
\label{eq<,2<p<3} \E|X+S_{\mathbf{X}}|^p\le A+\E\bigl|X+B^{1/2}Z\bigr|^p.
\end{equation}
Suppose now that $\sigma\ne0$.
Define the function $f$ by the formula $f(x):=\frac{\E|X/\sigma
+x|^p}{p(p-1)}$ for all $x\in\mathbb{R}$.
Using Lemma~\ref{lemunderint,simpl}, it is easy to see that $f''(x)=\E
|X/\sigma+x|^{p-2}$ for all $x\in\mathbb{R}$, and hence the function
$f$ is in the class $\mathscr{F}_p$ defined on page~515 from \cite{tyurinSPL}.
It follows\vspace*{1pt} by Theorem~2 in \cite{tyurinSPL}, and Jensen's inequality that
$\E|X+S_{\mathbf{X}}|^p\le\E|X+\sigma Z|^p+A\le\E|X+\sigma
Z+\sqrt{B-\sigma
^2}Z_1|^p+A=\E|X+B^{1/2}Z|^p+A$, where $Z_1\sim N(0,1)$.
So, inequality (\ref{eq<,2<p<3}) holds as well in the case $\sigma\ne0$.
Thus, the first supremum in (\ref{eq2<p<3}) is no greater than $A+\E
|X+B^{1/2}Z|^p$.

It remains to show that the second supremum in (\ref{eq2<p<3}) is no
less than $A+\E|X+B^{1/2}Z|^p$.
Recall (\ref{eqQ=}) and take any quadruple $(c_1,c_2,\lambda
_1,\lambda_2)\in Q_{p;A,B}$.
By Lemma~\ref{lemlim} with $G=\lambda_1\delta_{c_1}+\lambda_2\delta
_{c_2}$, one has a sequence $(\ZZ_n)$ in $\mathscr{X}_{p;A,B}$ such
that $S_{\ZZ_n}\stackrel{\mathrm{D}}\longrightarrow c_1\tPi
_{\lambda_1}+c_2\tPi_{\lambda_2}$.
By the Fatou lemma (Theorem~5.3 in \cite{billingsley}),
$
\liminf_n\E|X+S_{\ZZ_n}|^p\ge\E|X+c_1\tPi_{\lambda_1}+c_2\tPi
_{\lambda_2}|^p$, 
so that the second supremum in (\ref{eq2<p<3}) is no less than $\E
|X+c_1\tPi_{\lambda_1}+c_2\tPi_{\lambda_2}|^p$, for any
$(c_1,c_2,\lambda_1,\lambda_2)\in Q_{p;A,B}$.
So, by Proposition~\ref{propA+E} (whose proof does not rely on
Theorem~\ref{th2<p<3}), this supremum
is indeed no less than $A+\E|X+B^{1/2}Z|^p$.
\end{pf*}

\begin{pf*}{Proof of Proposition~\ref{propA+E}}
Let the quadruple $(c_1,c_2,\lambda_1,\lambda_2)\in Q_{p;A,B}$ vary
as in (\ref{eq=lim}), so that $c_1\to0$ and $|c_2|\to\infty$.
For $j\in\{1,2\}$, let $w_j:=c_j^2\lambda_j$, so that $w_1+w_2=B$,
$|c_1|^{p-2}w_1+|c_2|^{p-2}w_2=A$ and $c_1\tPi_{\lambda_1}+c_2\tPi
_{\lambda_2}\stackrel{\mathrm{D}}=Y_H$ with $H:=w_1\delta
_{c_1}+w_2\delta_{c_2}$.
It follows that $|c_1|^{p-2}w_1\le|c_1|^{p-2}B\to0$ and hence
$|c_2|^{p-2}w_2\to A$.
It remains to refer to
Proposition~\ref{proplim} (with $q=p$),
since $Y_{B\delta_0}\stackrel{\mathrm{D}}=B^{1/2}Z$.
\end{pf*}

\begin{pf*}{Proof of Corollary~\ref{corA+E}}
The first equality in (\ref{equnified}) follows immediately by
Theorems~\ref{th} and \ref{th2<p<3}.
Also, by Lemma~\ref{lemlim} with $G=\lambda_1\delta_{c_1}+\lambda
_2\delta_{c_2}$ and the Fatou lemma (Theorem~5.3 in \cite{billingsley}),
$\E|X+c_1\tPi_{\lambda_1}+c_2\tPi_{\lambda_2}|^p$ is no greater
than the second supremum in (\ref{equnified}), for each
$(c_1,c_2,\lambda_1,\lambda_2)\in Q_{p;A,B}$.
So, the last supremum in~(\ref{equnified}) is no greater than the
first two ones there.

On the other hand, the last supremum in (\ref{equnified}) is obviously
no less than the maximum in (\ref{eq}), and, by Proposition~\ref
{propA+E}, this supremum is no less than $A+\E|X+B^{1/2}Z|^p$.
So, by Theorems~\ref{th} and \ref{th2<p<3}, the last supremum in
(\ref{equnified}) is no less than the first two ones there.
\end{pf*}

\begin{pf*}{Proof of Theorem~\ref{thsymm}}
This proof is analogous to that of Theorem~\ref{th} and even
significantly simpler overall, since analogues of Propositions~\ref
{propcard=1} and~\ref{proplim} are not needed here.
In the proofs of the analogues of Propositions~\ref{propincr}, \ref
{prop<2}, and~\ref{propnogaus}, one should use the symmetrized
real-valued measure $\Delta(\dd u)+\Delta(-\dd u)$ in place of
$\Delta(\dd u)$.
\end{pf*}

\section*{Acknowledgement}
I am pleased to thank the referees for useful and stimulating comments.




%

\printaddresses
\end{document}